\documentclass[12pt]{amsart}
\usepackage{amssymb}
\usepackage{latexsym}
\usepackage{graphicx}
\usepackage{subfigure}
\usepackage{tikz}
\usepackage[linktocpage=true]{hyperref}
\usepackage{ulem}
\usepackage[left=3cm, right=3cm, top=3cm, bottom=3cm]{geometry}
\usepackage{color}

\newcommand{\re}[1]{(\ref{#1})}

\newcommand{\rl}[1]{Lemma~\ref{#1}}

\newcommand{\rp}[1]{Proposition~\ref{#1}}
\newcommand{\rt}[1]{Theorem~\ref{#1}}
\newcommand{\rd}[1]{Definition~\ref{#1}}
\newcommand{\rrem}[1]{Remark~\ref{#1}}

\hypersetup{ colorlinks   = true, 
urlcolor  = blue, 
linkcolor    = blue, 
citecolor   = blue 
}

\def\I{{\mathrm I}}

\def\om{{\omega}}

\def\lam{{\lambda}}
\def\del{{\delta}}

\def\al{{\alpha}}
\def\Ga{{\Gamma}}
\def\De{{\Delta}}
\def\ti{\tilde}
\def\ka{{\kappa}}

\def\beq{\begin{equation}}
\def\eeq{\end{equation}}

\def\Exp{\mathrm{Exp}}
\def\IM{\mathrm{Im}}
\def\Ker{\mathrm{Ker}}

\def\det{\mathrm{det}\ }

\def\id{\mathrm{Id}}
\def\bm{\begin{matrix}}
\def\em{\end{matrix}}

\newcommand{\V}{{\mathbb V}}

\newcommand{\R}{{\mathbb R}}

\newcommand{\C}{{\mathbb C}}

\newcommand{\N}{{\mathbb N}}
\newcommand{\PP}{{\mathbb P}}

\newcommand{\CD}{{\mathcal D}}
\newcommand{\CE}{{\mathcal E}}
\newcommand{\CF}{{\mathcal F}}
\newcommand{\CG}{{\mathcal G}}

\newcommand{\CH}{{\mathcal H}}
\newcommand{\CI}{{\mathcal I}}
\newcommand{\CL}{{\mathcal L}}

\newcommand{\CN}{{\mathcal N}}

\newcommand{\CW}{{\mathcal W}}
\newcommand{\CP}{{\mathcal P}}
\newcommand{\CS}{{\mathcal S}}

\newcommand{\CZ}{{\mathcal Z}}
\newcommand{\CU}{{\mathcal U}}
\newcommand{\CK}{{\mathcal K}}

\newcommand{\CQ}{{\mathcal Q}}
\newcommand{\CJ}{{\mathcal J}}
\newcommand{\CR}{{\mathcal R}}

\newtheorem{thm}{Theorem}[section]
\newtheorem{cor}[thm]{Corollary}
\newtheorem{prop}[thm]{Proposition}
\newtheorem{lemma}[thm]{Lemma}
\newtheorem{remark}[thm]{Remark}

\newtheorem{defi}[thm]{Definition}

\newcommand{\la}{\langle}
\newcommand{\ra}{\rangle}

\def\blu{  \color{blue}}
\definecolor{gainsboro}{rgb}{0.86, 0.86, 0.86}
 \definecolor{gray-light}{rgb}{0.75, 0.75, 0.75}
  	\definecolor{gray}{rgb}{0.5, 0.5, 0.5}

\newcommand{\RE}{\operatorname{Re}}

\begin{document}

\title[]{Local rigidity of actions of isometries on compact real analytic Riemannian manifolds}

\author{Laurent Stolovitch}
\address{Universit\'e C\^ote d'Azur, CNRS, Laboratoire J. A. Dieudonn\'{e}, 06108 Nice, France}
\email{Laurent.stolovitch@univ-cotedazur.fr}
\thanks{}

\author{Zhiyan Zhao}
\address{Universit\'e C\^ote d'Azur, CNRS, Laboratoire J. A. Dieudonn\'{e}, 06108 Nice, France}
\email{zhiyan.zhao@univ-cotedazur.fr}
\thanks{}

\begin{abstract}
In this article, we consider analytic perturbations of isometries of an analytic Riemannian manifold $M$. We prove that, under some conditions, a finitely presented group of such small enough perturbations is analytically conjugate on $M$ to the same group of isometry it is a perturbation of. Our result relies on a ``Diophantine-like" condition, relating the actions of the isometry group and the eigenvalues of the Laplace-Beltrami operator. Our result generalizes Arnold-Herman's theorem about diffeomorphisms of the circle that are small perturbations of rotations.
\end{abstract}

%



\maketitle

\section{Introduction and main result}

The aim of this article is to study the analytic rigidity of a group action by isometries on a compact real analytic Riemannian manifold $M$ (supposed to be connected and without boundary). Both $M$ and its Riemannian metric $g$ are supposed to be analytic. We consider a finitely presented group $G$, together with a group action by isometries $\pi$ on $M$, and an analytic group action by diffeomorphisms $\pi_0$ on $M$ which is a small perturbation of $\pi$. Our aim is to give conditions for $\pi_0$ and $\pi$ such that $\pi_0$ is analytically conjugated to $\pi$.

This problem takes its roots in the seminal articles of Arnold \cite{arnold-cercle}, Hermann \cite{hermann-ihes} and Yoccoz \cite{yoccoz-cercle1998} dedicated to analytic circle diffeomorphisms.
It was proved that, if such a diffeomorphism $F$ is a {\it small} perturbation of a rotation  $R_{\alpha}$ of a {\it Diophantine} angle $\alpha$ and if the rotation number of $F$ is also $\alpha$, then $F$ is analytically conjugated to $R_\alpha$. A similar statement was obtained in the smooth category by Moser \cite{Moser-cercle} for abelian groups of smooth circle diffeomorphisms. Global results (i.e., without the smallness assumption on the perturbation) are due to Herman \cite{hermann-ihes}, Yoccoz \cite{yoccoz-cercle} and to Fayad-Khanin \cite{FK-cercle} for single circle diffeomorphisms and abelian group of the latter, respectively. These ridigity  problems have a long history. See recent works in somehow different contexts in \cite{Khanin-ICM18,FK-ICM18,DK-duke}

In the present work, we consider a general compact real analytic manifold $M$ - that plays the role of the circle - on which acts a group of isometries, defined by a finite number of letters and relations. The isometries play the role of rotations. We define an appropriate notion of ``Diophantiness" for the isometries. This depends heavily on the geometry and the metric of the manifold $M$, as the spectrum of the Laplace-Beltrami operator is involved.
The appropriate notion that replaces the condition on ``rotation number of the perturbation" can be rephrased as ``the perturbation can be conjugated as close as one wishes to the unperturbed isometries". In different context, one would say ``the perturbation is formally conjugated or almost reducible to the unperturbed one". The point is then to prove that one can effectively achieve a genuine analytic conjugacy between the unperturbed and the small enough perturbed actions. 

Following notions in Section \ref{subsec_GroupAction}, we have that an action of a (finitely presented) group $G$ by isometries on a manifold $M$ induces an action on the tangent bundle $TM$. It gives rise to the Hochschild complex of cochains of $L^2$ vector fields~:
\begin{equation}\label{complex}
L^2(M,TM)\stackrel{d_0}{\longrightarrow}C^1(G, L^2(M,TM))\stackrel{d_1}{\longrightarrow}C^2(G, L^2(M,TM))\longrightarrow \cdots .
\end{equation}
We then introduce the self-adjoint box operator, which is fundamental for our purpose and is precisely defined in (\ref{def-box}), 
$$
\square:=d_0\circ d_0^*+d_1^*\circ d_1,
$$
with the adjoint being defined upon the $L^2-$scalar product on $M$.  Relating the spectral properties of $\square$ to that of the Laplace-Beltrami operator on the tangent bundle $\Delta_{TM}$, the ``Diophantine condition" is defined in Definition \ref{defi-dio}, while the definition of ``formal conjugacy" is given in Definition \ref{defi_formal_conj}.

We give below the main result of this article, the precise formulation of which is given in \rt{thm_conj}.

\begin{thm}\label{thmmain}
  Let $M$ be a real analytic compact manifold with an analytic Riemannian metric. Let $G$ be a finitely presented group and let $\pi$ be a Diophantine $G$-group action by analytic isometries on $M$. We assume that $\dim\Ker\;  \square <+\infty$. Let $\pi_0$ be an analytic $G$-group action by diffeomorphisms on $M$ which is sufficiently close to $\pi$. If $\pi_0$ is formally conjugated to $\pi$, then $\pi_0$ is analytically conjugated to $\pi$.
\end{thm}

\medskip

\noindent{\it Idea of proof.} In the circle diffeomorphism case, one first realizes the circle by the real line and the diffeomorphism by a Fourier series that can be extended holomorphically in a strip around the real axis in the complex plane. Then one conjugates the perturbation of the rotation on the strip to a much smaller perturbation of the rotation but on a narrower strip. The process can be continued each time through the limit, that is the rotation on some strip around $\R$.

We proceed roughly in the same way. We consider complex neighborhoods $M_{r}$, with $r>0$ sufficiently small, of $M=:M_0$, known as {\it Grauert tubes}. They were studied in a series of seminal articles by Sz\"oke \cite{Szo91,Szo95}, Lempert-Sz\"oke \cite{LeSz91}, Guillemin-Stenzel \cite{GS91,GS92} by considering a complex structure defined on the tangent bundle. On the other hand, Boutet de Monvel \cite{BdM78,Leb18,Zel12} defined them as a domain of holomorphic extension of the eigenvectors of the Laplace-Betrami operator of $M$. It happens that these notions coincide \cite{GLS96}.

These eigenvectors, and more precisely their counterpart in the space of analytic sections of the tangent bundle of $M$, together with their associated eigenvalues, play an important role in the present work. Indeed, we shall view every single perturbation $\pi_0(\gamma)$, with $\gamma$ the generator of $G$, of the isometry $\pi(\gamma)$ as $\Exp\{P_0(\gamma)\}\circ\pi(\gamma)$, where $P_0(\gamma)$ is an analytic vector field on $M$ and $\Exp$ denotes the exponential relative to the Riemannian connection of $M$.
We shall decompose $P_0(\gamma)$ along these eigenvectors, providing a kind of ``Fourier-like decomposition".
Such analytic vector field $P_0(\gamma)$ extends to a holomorphic vector field on some Grauert tube $M_{r_0}$.
Using the work of Boutet de Monvel, we shall only consider vector fields belonging to the {\it Hardy space of $M_{r_0}$} (see \rd{def-Hardy}). 
These are holomorphic vector fields over $M_{r_0}$ with $L^2$ boundary values, characterized by the decay of the coefficients of its ``Fourier-like decomposition" and depending on the ``radius" $r_0$ of the Grauert tube. This defines naturally a weighted $L^2-$norm \re{norme-ponderee}.

The group action $\pi$ by isometries gives rise to a self-adjoint operator $\square$, called ``box operator" and defined in \re{def-box}, on $(L^2(M,TM))^k$, with $k$ the number of generators and their inverses of $G$. The space $L^2(M,TM)$ can be decomposed into direct sum of irreducible finite dimensional subspaces $V_i$, each of which is contained in an eigenspace of the Laplace-Beltrami operator $\Delta_{TM}$ associated to its eigenvalue {\blu $-\lambda_i^2$}.   We shall then consider the ``linearized conjugacy equation" decomposing the vector field along these irreducible spaces. The ``Diophantine-like" condition \re{diophantine} given in Definition \ref{defi-dio} means that the eigenvalues of the box operator restricted to $V^k_i$ accumulate the origin not faster than a fixed power of the inverse of $\lambda_i$. 
For the group action $\pi$, it will be shown that, if the set of generators of $G$ is Diophantine in the sense of Dolgopyat \cite{dolgopyat} (see \rd{dioph-dolg}), then $\pi$ satisfies our Diophantine-like condition \re{diophantine}.

Essentially, the proof proceeds through an iterative scheme.  More precisely, for an $\varepsilon_m-$close perturbation of $\pi$ in a Grauert tube of ``radius" $r_m>\frac{r_0}2$ at the $m-\rm th$ iteration step, with $\varepsilon_m$ and $r_m$ sufficiently small, the ``Diophantine-like" property \re{diophantine} allows us to solve the linearized equation associated to the rigidity problem up to $\varepsilon_{m+1}$, which is an almost square of $\varepsilon_{m}$, size error plus a ``harmonic component", a priori of size $\varepsilon_m$. 
Moreover, the ``formal rigidity" assumption allows us to show that this harmonic component is indeed also of size $\varepsilon_{m+1}$. This allows us to conjugate the perturbation of the isometry to another one, the size of which is almost the square of the previous one, but on a narrower Grauert tube of ``radius" $r_{m+1}$ with $\frac{r_0}2<r_{m+1}< r_m$. See \re{seqs} for the precise choice of sequences of quantities in the iterative scheme.

\medskip

As a particularly interesting example of local analytic rigidity, we have the following theorem, which can be seen as an analytic version of Fisher's ``local rigidity" result \cite{Fisher}[Theorem 1.1] of Diophantine $G-$action by analytic isometries.
\begin{thm} \label{thm-geom0}
Let $\pi$ be a Diophantine $G-$action by analytic isometries on $M$ as above. Assume that the first cohomology group $H^{1}(G,L^2(M,TM)):=\Ker\, d_1/\IM\, d_0$ of the complex \re{complex} vanishes. Then any small enough analytic perturbation $\pi_0$ of $\pi$ is analytically conjugate to $\pi$.
\end{thm}
\noindent The definition of first cohomology will be given in Section \ref{subsec_GroupAction}.

\medskip

The following theorem can be seen as an analytic version of results by Moser \cite{Moser-cercle}, Karaliolios \cite{nikos-tore} and Petkovic \cite{petko-tore} relative to simultaneous conjugacy of a commutative family of perturbations of rotations on the torus to rotations. Let $\CG=\{e_1,\cdots, e_m\}$ be the canonical basis of $\mathbb{Z}^m$. 
\begin{thm}\label{thm-tore0}
Let $\pi$ be Diophantine $\mathbb{Z}^m$-action by rotations on the torus $\mathbb{T}^d$ : Let $\al_i\in \mathbb{R}^d$ be the rotation vector of the rotation $\pi(e_i)$. Assume there exist $c ,\tau>0$, such that for all $({\bf k},l)\in \mathbb{Z}^d\times \mathbb{Z}\setminus\{0\}$,$$\max_{1\leq i\leq m}|\la {\bf k},\al_i\ra-l|\geq \frac{c}{|{\bf k}|^{\tau}}.$$
Then any small enough analytic perturbation $\pi_0$ (isotopic to $\rm Id$) of $\pi$ such that, for each each $i$, the rotation vector $\al_i$ belongs to the convex hull of rotation set of $\pi(e_i)$, is analytically conjugate to $\pi$.
\end{thm}
Our results also applies to other examples, see Section \ref{exemple-sect}.

\medskip

\noindent{\it Description of the remaining of paper.} 
The remaining of paper is organized as follows. Section \ref{section-geom} is devoted to the geometric setting, definitions and properties of norms as well as the holomorphic counterpart of a lemma of Moser \cite{Mos69}.
Section \ref{subsec_GroupAction} is devoted to the setting of group actions, and gives the precise formulation of the main result.
Section \ref{section-box} contains properties of the box operator, and gives the estimate of solutions of the cohomological equations. Section \ref{sec-KAM} is devoted to a Newton scheme that proves the analyticity of conjugacy to the unperturbed group action.

\medskip

\noindent{\it Acknowledgment.} This work was stimulated by a work of David Fisher \cite{Fisher} with whom the first author had discussions around 2007 about it. Although not published, this article contains lot of interesting examples, in the smooth category. The first author thanks Charlie Epstein for having pointed out Boutel de Monvel's theory of holomorphic extension of eigenvectors, L\'{a}zl\'{o} Lempert, Matthew Stenzel and Robert Sz\"{o}ke for exchanges about Grauert tubes and also Bassam Fayad and Jonhattan DeWitt for exchanges on group actions. We thank David Fisher for having pointed out a mistake in Examples section \ref{exemple-sect} of our first version and for exchanges that followed.

\section{Real analytic Riemannian manifold and Grauert tube}\label{section-geom}

Let $M$ be a compact real analytic Riemannian manifold of dimension $n\geq 1$. The Riemannian metric, which is supposed to be real analytic, is defined by means of scalar product $\langle\cdot,\cdot\rangle_m$ on the tangent space $T_mM$ for every $m\in M$. We shall write, in local coordinates $(x_1,\cdots,x_n)$ over which the tangent bundle is trivialized,
$$\langle v,w\rangle_m=\sum_{1\leq i,j\leq n}g_{i,j}(x(m))v^i w^j,\qquad {\rm for} \  \  v=\sum_{i=1}^nv^i\frac{\partial}{\partial x_{i}}, \quad  w=\sum_{i=1}^nw^i\frac{\partial}{\partial x_{i}},$$
with the matrix $g(m)=(g_{i,j}(x(m)))_{1\leq i,j\leq n}$ positive definite at every $m\in M$,
 if $(m,v)$ and $(m,w)$ belong to $T_mM$.
Let $\al_m : T_mM\rightarrow T_m^*M$ be the isomorphism
$$\al_m(v)w:=\langle v,w\rangle_m,\qquad {\rm for} \  \  v,w\in T_mM.$$
This defines a scalar product on the cotangent bundle $T^*_mM$~:
$$\langle v^*,w^*\rangle_m:=\langle v,w\rangle_m, \qquad {\rm for} \  \ v^*=\al_m(v),\ w^*=\al_m(w).$$
The above is extended to an isomorphism
$\al_m: \wedge^pT_mM\rightarrow \wedge^pT^*_mM$ by
$$\alpha_m(v_1\wedge \cdots\wedge v_p):=\alpha_m(v_1)\wedge \cdots\wedge \alpha_m(v_p). $$ The scalar product induced on $\wedge^pT_mM$ is defined to be
$$\langle v_1\wedge \cdots\wedge v_p,w_1\wedge \cdots\wedge w_p\rangle_m:=\det (\langle v_i,w_j\rangle_m)_{1\leq i,j\leq p}.$$

\subsection{Exponential map and Grauert tube}

We recall (without proofs) some useful facts for real analytic Riemannian manifolds stated in \cite{GLS96}[Section 1].
First of all, according to Bruhat-Whitney theorem \cite{WB59} (see \cite{GLS96}[Lemma 1.2]), $M$ can be identified with a totally real submanifold of a complex analytic manifold $\tilde M$ of (real) dimension $2n$~: for all $m\in M$, there exists an open neighborhood $W$ of $m$ in $\tilde M$ and a holomorphic coordinate system $(z_1,\cdots,z_n)$ on $W$ such that
\beq\label{im}
W\cap M=\{q\in W : \IM z_1(q)=\cdots=\IM z_n(q)=0 \}.
\eeq
We recall a well known fact (see \cite{GLS96}[Corollary 1.3]).
\begin{prop}\label{holom-ext}\cite{Tom66} Let $M\hookrightarrow \tilde M$ be a totally real submanifold of a complex manifold $\tilde M$. Let $M'$ be a complex manifold and let  $f:M\rightarrow M'$ be a real analytic mapping. Then, there exists an open connected neighborhood $W$ of $M$ in $\tilde M$ and a unique holomorphic mapping $f^+:W\rightarrow M'$ such that $\left. f^+\right|_{M}=f$.
\end{prop}

For $m\in M$, let $B_m(0,r)\subset T_mM$ be the ball in $T_mM$ centered at $0$ and of radius $r$. Let $\Exp_m$ denotes the {\it exponential map} defined upon the Riemannian connection \cite{Hel01}[Chap. I, Section 6]. There exists $r(m)>0$ such that the mapping
$$
\Exp_m : B_m(0,r(m))\subset T_mM \rightarrow M
$$
is an analytic diffeomorphism onto its image. Moreover, the mapping $m\mapsto r(m)$ can be chosen lower semicontinuous.


Following \cite{GLS96}[Corollary 1.3], for $m\in M$, there exists an open connected neighborhood $W_m\subset T_mM\otimes\Bbb C$ and a unique holomorphic extension of $\Exp_m$ on $W_m$, still denoted by $\Exp_m$, to $\tilde M$.
Moreover, according to \cite{GLS96}[Theorem 1.5], there exists $0<r_*\leq\inf_{m\in M}r(m)$ such that for every $0<r<r_*$, the map
\beq\label{complex-exp}
\Phi : T^{r}M \rightarrow \tilde M,\quad \Phi(m,\xi)=\Exp_m\{{\rm i}\xi\}
\eeq
is an analytic diffeomorphism onto its image, where $$T^{r}M:=\{(m,\xi)\in TM : \|\xi\|_{g(m)}<r\}.$$
According to \cite{Szo91}[Theorem 2.2], \cite{GLS96}[Proposition 1.7] and \cite{LeSz91}, for any $0<r<r_*$, $T^{r}M$ admits a unique complex structure for which the complexified exponential
$$
T^{r}M \ni(m,\xi)  \mapsto \Exp_m\{{\rm i}\xi\} \in \Phi(T^{r}M)=: M_{r}
$$
is a biholomorphism. We shall write $TM^{\mathbb{C}}:=TM\otimes_{\mathbb{R}}\mathbb{C}$.

According to \cite{Gra58}(see also \cite{GS91,GS92}[Introduction]), there exists a non-negative smooth strictly plurisubharmonic function
\beq\label{rho}
\rho:  M_{r_*}\rightarrow [0,r_*[  \  \text{ with } \rho^{-1}(0)=M \text{ and } M_{r}=\rho^{-1}([0,r[), \ 0<r<r_*.
\eeq
Moreover, there exists an anti-holomorphic involution $\sigma :  M_{r_*}\rightarrow  M_{r_*}$ whose fixed point set is $M$ and $\rho(\sigma(q))=\rho(q)$ for all $q \in  M_{r_*}$.

Since the metric $g$ on $M$ is real analytic, it turns out that such a $M_{r}$ can be defined by a unique {\it real analytic} strictly plurisubharmonic function $\rho$ such that the K\"ahler form
$$
\omega:=\frac{\rm i}{2}\partial\bar\partial \rho= \frac{\rm i}{2}\sum_{1\leq i,j\leq n}\frac{\partial^2\rho}{\partial z_i\partial \bar z_j}dz_i\wedge d\bar z_j
$$ defines a K\"ahler metric on $M_{r}$, $0<r<r_*$,  \beq\label{kappa}\kappa:= \sum_{1\leq i,j\leq n}\frac{\partial^2\rho}{\partial z_i\partial \bar z_j}dz_i\otimes d\bar z_j,\eeq
which extends the Riemannian metric $g$ on $M$ according to the following theorem.

\begin{thm}\cite{GS91}[P.562]\label{GS-thm}
There exists a neighborhood $U$ of $M$ in $\tilde M$ and a unique real analytic solution $\rho$ on $U\setminus M$ of the complex Monge-Ampère equation
\beq\label{MA}
\det\left(\frac{\partial^2 \sqrt{\rho}}{\partial z_i\partial\bar z_j}\right)=0
\eeq
such that the inclusion map $(M,g)\hookrightarrow (\tilde M, \kappa)$ is an isometric embedding.
\end{thm}
\noindent Hence, the boundary $\partial M_{r}$ of $M_{r}$ is a compact real analytic manifold and $\overline{M}_r:=\rho^{-1}([0,r])$ is a compact K\"ahler manifold. The complex neighborhood $M_r$ of $M$ is called a {\it Grauert's tube}. 


Let us first extend the exponential map w.r.t. the metric $g$ at $m\in M$,
$$ \Exp_m : B_m(0,r(m))\subset T_mM\rightarrow M$$
  to the one w.r.t. the metric $\kappa$ at $q\in M_{r}$,
$$\Exp_{q} : B_{q}(0,r(q))\subset T_{q}^{(1,0)} M_{r}\rightarrow M_{r}.$$
Let $X$ be a real analytic vector field on $M$. According to \rp{holom-ext}, it extends to a holomorphic vector field $X$ on an open connected neighborhood $\CU$ of $M$ in $ M_{r}$, still denoted $ M_{r}$.
That is, $X$ is a holomorphic section of $T^{(1,0)} M_{r}$ over $ M_{r}$.
First of all, according to \cite{GLS96}[Proposition 1.9, 1.13], if $r$ is small enough, the analytic Riemannian metric $g$ uniquely extends to a non-degenerate holomorphic section $g^+\in \Ga^\om\left( M_{r}, BS(T^{(1,0)} M_{r})\right)$, where $BS(T^{(1,0)} M_{r})$ denotes the bundle of symmectric bilinear forms on the holomorphic vector fields $T^{(1,0)} M_{r}$ over $M_{r}$, that is $g^+$ defines a {\it holomorphic Riemannian metric} \cite{Lebrun}.
For each $q\in M_{r}$, we define $\Exp_{q}$ on the ball $B_q(0,r(q))\subset T^{(1,0)}_{q}M_{r}$ with respect to its K\"ahler metric $\kappa$ as follow~:
given a coordinate chart $(U,x)=(U,x_1,\ldots, x_n)$ of $M$ trivializing $TM$, let $(W,z_1,\ldots, z_n)$ be a holomorphic chart of $M_{r}$ extending a chart $(U,x)$ of $M$ as in $(\ref{im})$, with $W\cap M=U$ and $x_i=\RE z_i$ trivializing $T^{(1,0)} M_{r}$. Let us write
$$g(x(m))=\sum_{1\leq i,j\leq n} g_{i,j}(x(m))dx_i\otimes dx_j,\quad g^+(z(q))=\sum_{1\leq i,j\leq n} g_{i,j}^+(z(q))dz_i\otimes dz_j,
$$
where the matrices $(g_{i,j}(x(m)))_{1\leq i,j\leq n}$, $(g_{i,j}^+(z(q)))_{1\leq i,j\leq n}$ are invertible for each point  $m\in M$ and $q\in M_{r}$ respectively. We recall that the geodesics on $M$ are solutions of the (real time) differential equation, in a coordinate chart~:
$$
\ddot x_j = \sum_{1\leq k,l\leq n}\Ga^j_{k,l}(x)\dot x_k\dot x_l,\quad j=1,\cdots, n,
$$
where, $\Ga^j_{k,l}$ denotes the Christoffel symbol defined by
$$
\Ga^j_{k,l}(x):=\frac{1}{2}\sum_{1\leq m\leq n} g^{j,m}(x) \left(\frac{\partial g_{m,k}}{\partial x_l}-\frac{\partial g_{k,l}}{\partial x_m}+\frac{\partial g_{l,m}}{\partial x_k}\right),
$$
and $(g^{j,m}(x))$ denotes the inverse matrix of $(g_{j,m}(x))$. 
Following \cite{Lebrun}[1.17, 1.18], let us consider the {\it holomorphic differential equation} with complex time:
\beq\label{geod-complex}
\ddot z_j = \sum_{1\leq k,l\leq n}\Ga^{+j}_{k,l}(z)\dot z_k\dot z_l,\quad j=1,\cdots, n,
\eeq
with $\Ga^{+j}_{k,l}$ defined as
$$
\Ga^{+j}_{k,l}(z):=\frac{1}{2}\sum_{1\leq m\leq n} (g^+)^{j,m}(z) \left(\frac{\partial g_{m,k}^+}{\partial z_l}-\frac{\partial g_{k,l}^+}{\partial z_m}+\frac{\partial g_{l,m}^+}{\partial z_k}\right),
$$
and $\left((g^+)^{j,m}(z)\right)$ the inverse matrix of $\left(g_{j,m}^+(z)\right)$. For any $q_0\in W$ and $(q_0,\xi)\in T_{q_0}^{(1,0)}M_{r}$, such that $(z_0,\xi)\in \Delta_1^n\times \mathbb{C}^n$ with $z_0=z(q_0)$, there exists a unique complex curve, a {\it complex geodesic}, $t\in D_{z_0,\xi}\mapsto z(t)=\Phi(t,z_0,\xi)$ with $(z(0),\dot z(0))=(z_0,\xi)$,
  solution of $(\ref{geod-complex})$. Here $D_{z_0,\xi}$ denotes a complex neighborhood of $0$ in $\Bbb C$ that depends on the point $(z_0,\xi)$. As in the real case, the form of Eq. \re{geod-complex} allows us to write $z(t)=\Psi(z_0,t\xi)$; it is holomorphic for $z_0\in\Delta_1^n$ and $t$ small complex number. Hence, $\Psi(z_0,\xi)$ is holomorphic for  $z_0\in\Delta_1^n$, and  $\xi$ in the complex ball 
in $\mathbb{C}^n$, centered at $0$ and of sufficiently small radius $\del$ w.r.t the K\"ahler metric $\kappa$~: $|\xi|_{\ka(q_0)}<\del$, $z(q_0)=z_0$.
Furthermore, it satisfies
\beq\label{basic-flow}
\Psi(z_0,0)=z_0,\quad D_{\xi}\Psi(z_0,0)=\id.
\eeq
Hence, for some holomorphic map $\varphi(z_0,\xi)$ satisfying $D_{\xi}\varphi(z_0,0)=0$, we have
\beq\label{exp-def}
z(t)= \Psi(z_0,t\xi)=z_0+t\xi+\varphi(z_0,t \xi).
\eeq
Taking a finite covering of $M_r$ by open sets, there exists an $a>0$ such that the solution $z(t)=\Phi(t,z(q),\xi)=\Psi(z(q),t\xi)$ is holomorphic for $|t|< 2$, $q\in  M_{r}$ and $|\xi|_{\kappa(q)}<a$. 
Let $(q,\xi)\in  T^{(1,0)}_{q} M_{r}$, be such a point (i.e. $|\xi|_{\ka}< a$) with $q\in W$. We define the {\it complex exponential map} $\Exp_{q}\{\xi\}$ to be the time$-1$ of this complex flow. It is the point of $M_{r}$ whose expression in the coordinate chart $W$ is
\beq\label{complex-exp}
\Exp_{q}\{\xi\}:= \Psi(z(q),\xi).
\eeq

Let $q\in W$ of sufficiently small coordinate $z(q)$. Let $(\xi,\eta)\in \mathbb{C}^{n}\times \mathbb{C}^n$ be small enough so that, 
$\Psi(\Psi(z(q),\xi),\eta)$ is well defined. According to \re{basic-flow}, there is a unique holomorphic map $(z,\xi,\eta)\mapsto P(z,\xi,\eta)=:\zeta\in \mathbb{C}^n$ 
that solves the equation  $\Psi(\Psi(z,\xi),\eta)=:\Psi(z,\zeta)$ for $(\xi,\eta)$ in small neighborhood of $0$ in $\mathbb{C}^{2n}$ and $z$ in a neighborhood of $z(q)$. 
 Furthermore, there is a holomorphic map $\theta(z,\xi,\eta)\in \mathbb{C}^n$ such that
$$ \zeta =\xi+\eta+\theta(z,\xi,\eta),\quad \theta(z,0,\eta)=\theta(z,\xi,0)=0,\quad \|D_\xi\theta(z,\xi,\eta)\|_{\ka}\leq c|\eta|_{\ka}.$$

\begin{remark}\label{complex-riemann}
We stress that, since $M_r$ is a K\"alher manifold, the unique Hermitian connection which is compatible with the metric and the complex structure on $T^{(1.0)}M_r$  coincide with the Riemannian connection of $\kappa$ \cite{zheng-book}[Propostion 7.9, Definition 7.13].
	\end{remark}

\subsection{Moser's lemma for Riemannian geometry}\label{sec_Moser}
The following proposition is an adaptation to our holomorphic context of those in the article of Moser \cite{Mos69}. Their proofs are identical: we just consider the holomorphic extension \re{complex-exp} of the ``Riemannian exponential", together with the holomorphic implicit theorem instead of the smooth implicit theorem.

Given $0<r<r_*$, let the set of holomorphic sections of $T^{(1,0)} M_{r}$ over $M_{r}$, that is holomorphic vector fields, be denoted by $\Gamma_{r} = \Gamma( M_{r},T^{(1,0)} M_{r})$,
equipped with the norm
$$
 |v|_{0,r}:=\sup_{q\in M_{r}}|v(q)|_{\ka},\quad v\in \Gamma_{r}.
$$
There is an analytic trivializing atlas with a finite covering patches $\left\{U_i,x^{(i)}\right\}_i$ of $M$ that extends to a holomorphic atlas $\left\{W_i,z^{(i)}\right\}_i$ of $M_{r_*}$ as in \re{im} and such that $z^{(i)}\in \Delta_1^n$ (In what follows, $z^{(i)}$ stands for $z^{(i)}(q)$ with $q\in W_i$). For $v\in\Gamma_{r} $, restricting to $W$, one of these coordinates patches on which $z\in \Delta_1^n$ and writing $\tilde v(z)=\sum_{1\leq j\leq n} \tilde v_j(z)\frac{\partial}{\partial z_j}$ the expression of $v$ in this coordinate patch, we set
\begin{equation}\label{norm01r}
 |v|_{0,r}:=\max_i\sup_{q\in W_i\cap M_{r}}|\tilde v(z^{(i)}(q))|_{\kappa}, \quad |v|_{1,r}:=\max_i\sup_{q\in W_i\cap M_{r}}\sup_{\substack{\zeta\in \mathbb{C}^n,\\|\zeta|\leq 1}}|D\tilde v(z^{(i)}(q))\zeta|_{\kappa},
\end{equation}
and $\|v\|_{1,r}:=|v|_{0,r}+|v|_{1,r}$. Moreover, we define the norms on the set of analytic sections of $TM$ on $M$, analytic vector fields on $M$, denoted by $\Gamma^{\omega}:=\Gamma^{\omega}(M,TM)$, as in \cite{Mos69}, \beq\label{norm-moser}
|v|_{0,M},\quad |v|_{1,M}\quad \|v\|_{1,M}:=|v|_{0,M}+|v|_{1,M}, \quad v\in \Gamma^{\omega},
\eeq in the similar sense as (\ref{norm01r}). In what follows, the $|\cdot|_{0,r'}-$norm with $r'=0$ means the $|\cdot|_{0,M}-$norm.
Since every $v\in \Gamma^{\omega}$ can be holomorphically extended to $M_r$ for some $0<r<r_*$, let $\Gamma_r^{\omega}\subset \Gamma_r$ be the set of holomorphic extensions to $M_r$ of elements in $\Gamma^{\omega}$.
It is obvious that $|v|_{0,M}\leq |v|_{0,r}$ for $v\in \Gamma^\om_r$.

In the following, the inequality with ``$\lesssim$" means boundedness from above by a positive constant depending only on the manifold $M_{r_*}$ and the K\"ahler metric $\ka$ on $M_{r_*}$ but independent of other factors.

\begin{prop}\label{lemMoser1}
The following assertions hold true.
\begin{enumerate}
\item [(i)] There exists sufficiently small $\varepsilon_*>0$, depending only on the manifold $(M_{r_*},\ka)$, such that for $v,w\in\Gamma_{r}$ with $0<r<r_*$, if $|v|_{0,r}$, $\|w\|_{1,r}<\varepsilon_*$, then there exists $s_1(w,v)\in \Gamma_{r'}$ for any $r'\in ]0,r[$ such that,
\begin{equation}\label{eq_lemMoser1}
\Exp \{w\}\circ \Exp \{v\} = \Exp \{w+v+s_1(w,v)\},
\end{equation}
satisfying $s_1(w,0)=s_1(0,v)=0$, and, for any $r'\in ]0,r[$, any 
$v_1,v_2\in\Gamma_{r}$ with $|v_1|_{0,r'}$, $|v_2|_{0,r'}<\varepsilon_*$,
$$|s_1(w,v_1)-s_1(w,v_2)|_{0,{r'}}\lesssim |w|_{1,r}|v_1-v_2|_{0,{r'}}.$$
\item [(ii)] There exists sufficiently small $\varepsilon_{M}>0$, depending only on the manifold $(M,g)$, such that for for $v,w\in\Ga^{\om}$ with $|v|_{0,M},  \|w\|_{1,M}<\varepsilon_{M}$, then there exists $s_1(w,v)\in \Gamma^{\om}$, such that (\ref{eq_lemMoser1}) is satisfied and
$$|s_1(w,v)|_{0,M}\lesssim \|w\|_{1,M} |v|_{0,M},\quad |s_1(w,v)|_{1,M}\lesssim \|w\|_{1,M} (1+|v|_{1,M}).$$

\end{enumerate}
\end{prop}

\begin{remark}\label{rem_s1}
In the above proposition, if $\tilde v, \tilde w\in \Gamma_r^\om$ are holomorphic extensions to $M_r$ of $v,w\in\Gamma^\om$, then $s_1(\tilde w,\tilde v)\in\Gamma_r^\om$ is the holomorphic extension of $s_1(w,v)\in\Gamma^\om$.
\end{remark}

Based on Lemma 1 in \cite{Mos69}, the assertion (i) is its holomorphic version on $M_r$, and the assertion (ii) gives the estimate on the derivatives of $s_1(w,v)$ on $M$. Both assertions can be deduced readily from the proof in \cite{Mos69}.
For completeness, we give a proof of this proposition in Appendix \ref{app_proof}.

\subsection{Spaces of sections of vector bundles}

Let $E$ be an analytic vector bundle over $M$. We shall denote by $\Ga^{\om}(M,E)$ (resp. $\Ga^{k}(M,E)$, $k\in \Bbb N\cup\{\infty\}$) be the space of analytic (resp. $C^k-$smooth) sections of $E$. If $E$ admits an analytic scalar product $\langle\cdot,\cdot\rangle_{E}$, then we define the scalar product on the space of section to be
$$
\langle v,w\rangle:=\int_M\langle v(x),w(x)\rangle_{E,x} \,  d{\rm vol}(x),
$$
where $d{\rm vol}$ is a volume element which can be expressed in local coordinates~:
$$
d{\rm vol}(x)= \sqrt{\det g_{i,j}(x)} \,  dx_1\cdots dx_n.
$$
Let $L^2(M,E)$ denotes the completion of $\Ga^{\infty}(M,E)$ with respect to this scalar product. It is the Hilbert space of $L^2$ sections of $E$.
The metric $g$ extends to a K\"ahler metric $\kappa$ on $M_{r}$.  Hence, for any holomorphic sections $v,w\in \Gamma_{r_*}$, $0<r<r_*$,
$$
\la v,w\ra:=\int_{M_{r}}\la v(z),w(z)\ra_{\ka} \, \frac{\om^n(z)}{n!}.
$$

From the de Rham complex, we construct a complex on the space of smooth sections of multi-vector fields as follow~:
$$
\begin{array}{cccccc}
	\Ga^{\infty}(M,\Bbb R) & \stackrel{d_0}{\longrightarrow}& \Ga^{\infty}(M,T^*M) & \stackrel{d_1}{\longrightarrow}& \Ga^{\infty}(M,\wedge^2 T^*M) & \stackrel{d_2}{\longrightarrow}\cdots\\
	\parallel & &\downarrow\alpha^{-1} & & \downarrow\alpha^{-1}& \\
	\Ga^{\infty}(M,\Bbb R) & \stackrel{\ti d_0}{\longrightarrow}& \Ga^{\infty}(M,TM) & \stackrel{\ti d_1}{\longrightarrow}& \Ga^{\infty}(M,\wedge^2 TM) & \stackrel{\ti d_2}{\longrightarrow}\cdots\\
\end{array}  \  \  \  \  .
$$
The first differentials are defined to be $\ti d_0 := \al^{-1}\circ d_0$ and $\ti d_1 := \al^{-1}\circ d_1\circ\alpha $. We shall call the $L^2$ extension of this complex the 
``tangential complex" of $M$.
Since the de Rham complex is elliptic (the complex of the associated symbols is exact), the same holds true for the tangential complex. We then define the Laplacian on the tangent bundle to be the self-adjoint operator
$$
\De_{TM} := \ti d_0\circ \ti d_0^*+ \ti d_1^*\circ \ti d_1.
$$
According to classical generalized Hodge theory for elliptic complex \cite{Wel80}[Section 5, Theorem 5.2], \cite{Dem96}[Theor\`eme 3.10, Corrolaire 3.16] or \cite{Ros97}, there exists an orthonormal basis $({\bf e}_j)_{j\geq 0}$ of eigenvectors of $-\De_{TM}$ in $L^2(M, TM)$, the associated eigenvalues are $(\tilde\lam_j^2)_{j\geq 0}$ with
\begin{equation}\label{tildelambda}
\tilde\lambda_0=0\leq\tilde\lam_1\leq\tilde\lam_2 \leq \cdots, \quad \lim_{j\to\infty}\tilde\lam_j=+\infty.\end{equation}
 Each eigenvalue is of finite multiplicity and $+\infty$ is the only accumulation point. Moreover, according to Weyl asymptotic estimate, we have that
\begin{equation}\label{asymp-lamdak}
\#\left\{j\in\N:\tilde\lambda^2_j \leq \lambda\right\}\sim a_0 \lambda^{\frac{n}{2}},\quad \tilde\lambda^2_j\sim b_0 \cdot j^{\frac2n}  \quad {\rm as}  \   \ j\to \infty
\end{equation}
for some constants $a_0, b_0>0$ depending only on the Riemannian manifold $(M,g)$.
(see e.g., \cite{B86}[Page 70]). The eigenvectors ${\bf e}_j$ are real analytic on $M$ \cite{treves-book-analytic}[Theorem 4.1.2].

\subsection{Analyticity of eigenvectors of Laplacian $\De_{TM}$}

We follow and recall the result of Boutet de Monvel \cite{BdM78}(see also \cite{GLS96}[Section 1 and 2] and \cite{Leb18}).

Let us consider the elliptic analytic pseudo-differential operator of order $1$, $|\De_{TM}|^{\frac12}$. Due to the classical elliptic theory, the eigenvectors of $|\De_{TM}|^{\frac12}$ (which are the same as those of $\De_{TM}$) are, in fact, real analytic. They can be considered as restrictions to $M$ of holomorphic sections on a same neighborhood of the $M$ in a complexified manifold $\tilde M$ of $M$.
\begin{defi}\label{def-Hardy}
Let the {\bf Hardy space} $\tilde H_r^2=\tilde H^2(M_{r},T^{(1,0)}M_{r})$ be the space of holomorphic sections of $T^{(1,0)}M_{r}$ over $M_{r}$ whose restriction (in the sense of distribution) to $\partial M_{r}$ belongs to $L^2(\partial M_{r},T^{(1,0)}M_{r})$, associated to
the {\bf Hardy product}
\beq\label{h2product}
\la f,h \ra_{\tilde H_r^2}:=
\int_{\partial M_{r}}\la f(q),h(q)\ra_{\ka}d\mu_{r}(q),\qquad f, h \in  \tilde H_r^2,
\eeq
and the {\bf Hardy norm}
\beq\label{h2norm}
	\|f\|_{\tilde H_r^2}:=\|f|_{\partial M_{r}}\|_{L^2(\partial M_{r})}=\left(\int_{\partial M_{r}}\la f(q),f(q)\ra_{\ka}d\mu_{r}(q)\right)^{\frac12},\qquad f\in\tilde H_r^2.
\eeq
Here, $d\mu_{r}$ denotes the ``surface measure" obtained by restriction of $\frac{\om^n(z)}{n!}$ to the real analytic level set $\rho=r$.
More generally, for $\nu\in\N^*$, the Hardy product and the Hardy norm on $(H_r^2)^\nu$ are defined as
\begin{eqnarray}
\la f,h \ra_{\tilde H_r^2}&:=&\sum_{1\leq l\leq \nu} \la f_l,h_l \ra_{\tilde H_r^2},\qquad f=(f_l)_{1\leq l\leq \nu}, \quad h=(h_l)_{1\leq l\leq \nu}\in (\tilde H_r^2)^\nu,\label{h2product-nu}\\
\|f\|_{\tilde H_r^2}^2&:=& \sum_{1\leq l\leq \nu}\|f_l\|_{\tilde H_r^2}^2, \qquad f=(f_l)_{1\leq l\leq \nu}\in (\tilde H_r^2)^\nu.\label{h2norm-nu}
\end{eqnarray}
Moreover, let the subspace
$\left(H_r^2\right)^\nu\subset \left(\tilde H_r^2\right)^\nu$ be
$$\left(H_r^2\right)^\nu:=\left\{f\in \left(\tilde H_r^2\right)^\nu: \left. f\right|_M \in (\Gamma^{\omega})^\nu=(\Gamma^{\omega}(M,TM))^\nu\right\},$$
equipped with the induced Hardy product.
\end{defi}

In what follows, we shall use the following ``vector-valued" version of Boutet de Monvel's theorem. It is obtained verbatim from its proof given by Stenzel \cite{stenzel-poisson}(or by Lebeau \cite{Leb18}) using the Heat kernel on sections of the tangent bundle (see \cite{gilkey}[Section 1.6.4, P.54]) instead of on functions. Indeed, recalling the definition of $\alpha:TM\rightarrow T^*M$ (see the beginning of Section \ref{section-geom}), the kernel
$$
K(t,x,y):=\sum_{k\geq 0}e^{-t\tilde\lambda_k}{\bf e_k}(x)\otimes \alpha({\bf e_k}(y))
$$
satisfies the elliptic system $(-2\partial_t^2+\I\otimes \De_{T^*M}+\Delta_{TM}\otimes \I)K=0$. Hence, it is analytic on $\mathbb{R}_+^*\times M\times M$ \cite{treves-book-analytic}[4.1.4].


\begin{thm}\label{theo-boutet}\cite{BdM78}
Let $f\in L^2(M,TM)$ be a global $L^2$ section of $TM$. It can be decomposed along the eigenvectors of $|\De_{TM}|^{\frac12}$~:
\beq\label{decomp_L2}
f=\sum_{j\geq 0} f_j{\bf e}_j, \quad f_j\in \R.\eeq Moreover, for $0<r<r_*$, $f$ extends to a section $\ti f\in H_r^2$ if and only if
	\beq\label{exp-decay-thm}
	\sum_{j\geq 0} \left|f_j\right|^2 e^{2 r\tilde\lam_j} (1+\tilde\lam_j)^{-\frac{n-1}{2}}<+\infty.
	\eeq
\end{thm}

\smallskip

In the sequel, the extension $\ti f\in H_r^2$ of $f\in L^2(M,TM)$ will be still denoted by $f$ since they are identified through the sequence of coefficients $(f_j)_{j\in\N}$ satisfying (\ref{exp-decay-thm}). Comparing with (\ref{tildelambda}), we reorder the distinct eigenvalues of $|\De_{TM}|^{\frac12}$ as
$$
\lambda_0=0<\lam_1<\lam_2 < \cdots, \quad \lim_{i\to\infty}\lam_i=+\infty.
$$
Taking account of the multiplicity of eigenvalues, we thus have the decomposition 
\begin{equation}\label{decomp_lambdas}
\N=\bigcup_{i\geq 0} I_i:= \bigcup_{i\geq 0}  \left\{j\in \N:\tilde\lam_j=\lambda_i\right\}
\end{equation}
Then (\ref{decomp_L2}) and (\ref{exp-decay-thm}) can be respectively rewritten as
\begin{equation}\label{form_multiplicity}
f=\sum_{i\geq 0}\sum_{j\in I_i} f_j {\bf e}_j,\quad \sum_{i\geq 0} e^{2 r\lam_i} (1+\lam_i)^{-\frac{n-1}{2}} \sum_{j\in I_i} \left|f_j\right|^2 <+\infty.
\end{equation}

 \begin{defi}\label{defi_norm_L2pon}
For $f\in L^2(M,TM)$, let the {\bf $L^2-$norm} be
\begin{equation}\label{norme-L2}
  \|f\|_{L^2}:=  \left(\int_{M}\la f(x),f(x)\ra_{g} d{\rm vol}(x)\right)^{\frac12}.
 \end{equation}
 For $f\in H_r^2$, $0<r<r_*$, let the {\bf weighted $L^2-$norm} be
\begin{equation}\label{norme-ponderee}
\|f\|_{r}:=\left(\sum_{j\geq 0}\left| f_j \right|^2 e^{2 r\tilde\lam_j} (1+\tilde\lam_j)^{-\frac{n-1}{2}}\right)^{\frac12},\quad f=\sum_{j\geq 0}f_j{\bf e}_j\in H^2_{r}.
\end{equation}
For the vector of sections, the norms are naturally defined as in \re{h2norm-nu}.
\end{defi}

\begin{remark}\label{rmk_L2} For $u=(u_l)_{1\leq l\leq\nu}\in (\Ga^\om)^\nu$, we have $u\in (L^2(M,TM))^\nu$ with
$$\|u\|_{L^2}=  \left(\sum_{1\leq l\leq \nu} \int_{M}\la u_l(x),u_l(x)\ra_{g}d{\rm vol}(x)\right)^{\frac12}\lesssim  \left(\sum_{1\leq l\leq \nu}|u_l|^2_{0,M}\right)^{\frac12}  =: |u|_{0,M}.$$
\end{remark}
In the following, let $\CI$ be a closed sub-interval of $]0,r_*[$.
The inequality with ``$\lesssim$" means boundedness from above by a positive constant uniform on $\CI$ (independent of the choice of $r\in \CI$) depending only on the manifold $(M_{r_*},\ka)$, and the inequality with ``$\simeq$" means such boundedness from above and below.

\begin{prop}\label{prop_estim-norm} The following assertions hold true for any $r\in {\CI}$.
\begin{enumerate}
\item [(i)] For every $v\in H_r^2$, $\|v\|_{r}\simeq \|v\|_{H^2_{r}} $, and for $u,v\in H_r^2$, $|\la u,v \ra_{H_r^2} |\lesssim \|u\|_r\cdot \|v\|_r$.

\smallskip

\item [(ii)] Given $v=\sum_{j\geq 0} v_j{\bf e}_j\in H^2_{r}$, the coefficients $\{v_j\}_{j\geq 0}$ satisfy that
\begin{equation}\label{decay_coeff}
|v_j|\leq \|v\|_{r}e^{- r\lam_i} (1+\lam_i)^{\frac{n-1}{4}},\qquad \forall \ j\in I_i.
\end{equation}
On the other hand, if the sequence $\{v_j\}_{j\geq 0}$ satisfies that
\begin{equation}\label{decay_coeff-inverse}
|v_j|\leq \CD e^{- r\lam_i} (1+\lam_i)^{\frac{n-1}{4}},\qquad \forall \ j\in I_i,
\end{equation}
for some constant $\CD>0$, then $v=\sum_{j\geq 0} v_j{\bf e}_j\in H^2_{r'}$ for any $r'\in ]0,r[$ with
\begin{equation}\label{esti-norm-vr'}
\|v\|_{r'}\lesssim \frac{\CD}{(r-r')^{\frac{n}{2}}}.
\end{equation}

\smallskip

\item [(iii)] For every $v\in H_r^2$, we have
\begin{eqnarray}
\|v\|_{r'}&\lesssim&|v|_{0,r},\qquad \forall \  r'\in ]0,r[. \label{supnormequiv-1}\\
|v|_{0,r'}&\lesssim&\frac{\|v\|_{r}}{(r-r')^{3n}},\qquad \forall \  r'\in [0,r[,\label{supnormequiv-2}
\end{eqnarray}
	
\smallskip

\item [(iv)] There are natural continuous embeddings $H^2_{r}\hookrightarrow H^2_{r'}$, $r'\in ]0,r[$.
	\end{enumerate}
\end{prop}

\noindent
{\it Proof of (i).}
In view of the inequality (1.7) in Theorem 1.1 of \cite{Leb18}, we have that, for any $ r\in ]0,r_*[$, there exists a constant $c_{r}>0$ such that
$$c_{r}^{-1}\|v\|_{H_r^2}\leq \|v\|_{r}\leq c_{r}\|v\|_{H_r^2},\qquad \forall \ v\in H_r^2.$$
Then, according to Proposition 5.4 of \cite{Leb18}, the constant $c_r$ in the above inequality can be chosen uniformly for $r$ in the sub-interval $\CI\subset ]0,r_*[$, which means the uniform equivalence between Hardy norm and weighted $L^2-$norm.
Through this equivalence and Definition \ref{def-Hardy}, we obtain that
\begin{eqnarray*}
|\la u,v \ra_{H_r^2} |&=&\left|\int_{\partial M_{r}}\la u(z),v(z)\ra_{\ka} \, d\mu_{r}(z)\right|\\
&\leq&\left|\int_{\partial M_{r}}\la u(z),u(z)\ra_{\ka} \, d\mu_{r}(z)\right|^{\frac12} \left|\int_{\partial M_{r}}\la v(z),v(z)\ra_{\ka} \, d\mu_{r}(z)\right|^{\frac12}\\
&\lesssim&  \|u\|_r \cdot \|v\|_r.
\end{eqnarray*}

\smallskip

\noindent
{\it Proof of (ii).} The estimate (\ref{decay_coeff}) follows immediately from Theorem \ref{theo-boutet} and its alternative form \re{form_multiplicity}. On the other hand, from the definition (\ref{norme-ponderee}) of weighted $L^2-$norm $\|\cdot\|_{r}$, and, under the assumption (\ref{decay_coeff-inverse}), $v=\sum_{j\geq 0} v_j{\bf e}_j$ satisfies that
$$\|v\|_{r'}^2  =   \sum_{j\geq 0}\left|v_j\right|^2 e^{2r'\tilde\lam_j} (1+\tilde\lam_j)^{-\frac{n-1}{2}}
		\leq  \CD^2 \sum_{j\geq 0}  e^{-2(r-r')\tilde\lam_j}.$$	
According to the asymptotic estimate (\ref{asymp-lamdak}), there is a constant $a$, depending only on the Riemannian manifold $(M,g)$ such that $\tilde\lam_j\geq  a j^{\frac{1}{n}}$. Hence, by successive integration by parts, we obtain (\ref{esti-norm-vr'}) through
$$\sum_{j\geq 0} e^{-2(r-r')\tilde\lam_j}\leq \sum_{j\geq 0} e^{-2a(r-r') j^{\frac{1}{n}}}\lesssim \int_0^{+\infty}  e^{-2a(r-r') t^{\frac{1}{n}}} dt \lesssim (r-r')^{-n}.$$

\smallskip

\noindent
{\it Proof of (iii).} By the definition of Hardy norm in \re{h2norm}, as well as the assertion (i), we have
$$\|v\|_{r'}\simeq \|v\|_{H^2_{r'}}\leq   \CS_{\CI}  \sup_{z\in \partial M_{r'}}|v(z)|_{\kappa}\leq   \CS_{\CI}|v|_{0,r},\quad \forall \  r'\in ]0,r[$$
where $\CS_{\CI}:=\sup_{r\in\CI} \CS_{r}$ with $\CS_{r}$ the ``surface size" of $\partial M_{r}$.

In view of Proposition 2.1 of \cite{GLS96}, we have, for every $r'\in[0,r_*[$,
$$
|{\bf e}_j|_{0,r'}\lesssim (1+\tilde\lam_j)^{n+1}e^{r'\tilde\lam_j},  \footnote{We also mention an improved estimate due to Zelditch \cite{Zel12}[Corollary 3] which allows to replace $(1+\tilde\lam_j)^{n+1}$ by $(1+\tilde\lam_j)^{\frac{n+1}{4}}$.}$$	
where the $|\cdot|_{0,r'}-$norm with $r'=0$ means the $|\cdot|_{0,M}-$norm defined in \re{norm-moser}.	
Therefore, by \re{decay_coeff}, for $v\in H_r^2$, for $r'\in [0,r[$,
$$|v|_{0,r'}  \leq \sum_{j\geq 0} |v_j||{\bf e}_i|_{0,r'}
\leq\|v\|_{r}  \sum_{j\geq0} e^{- (r-r')\tilde\lam_j} (1+\tilde\lam_j)^{\frac{5n+3}{4}}
\lesssim\frac{\|v\|_{r}}{(r-r')^{5n}}. $$
Indeed, in view of the asymptotic estimate (\ref{asymp-lamdak}), we have
\beq\label{sum}
\sum_{j\geq0}(1+\tilde\lam_j)^{\frac{5n+3}{4}} e^{-(r-r')\tilde\lam_j}\lesssim \sum_{j\geq1} j^{\frac{1}{n}\cdot \frac{5n+3}{4}} e^{-a(r-r') j^{\frac{1}{n}}},
\eeq
and, for the function $h_b:t\mapsto t^{\frac{c}{n}}e^{-bt^{\frac{1}{n}}}$, $b>0$ and $c=\frac{5n+3}{4}$, we have
$$\max_{t\in \R_+^*}h_b(t)=  h_b\left( \left(\frac{nc}{b}\right)^{n}\right) = \frac{d_n}{b^{\frac{5n}{4}+\frac{3}{4}}},$$
with a constant $d_n>0$ depending only on $n$. By successive integration by parts on $[0,+\infty[$, the sum \re{sum} is bounded as
$$
\sum_{j\geq1} j^{\frac{1}{n}\cdot \frac{5n+3}{4}} e^{-a(r-r') j^{\frac{1}{n}}}\lesssim \max_{t\in \R_+^*}h_{a(r-r')}(t) + \int_{0}^{+\infty} t^{\frac{5}{4}+\frac{3}{4n}} e^{-a(r-r')t^{\frac{1}{n}}} dt \lesssim \frac{1}{(r-r')^{3n}}.$$

\smallskip

\noindent
{\it Proof of (iv).}
By Definition \ref{def-Hardy} and Remark \ref{rmk_L2}, any $f=\sum_{j\geq 0} f_j{\bf e}_j\in H^2_{r}$, $ r\in \CI$, is an element of $L^2(M,TM)$ when it is restricted to $M$.
Moreover, for $r'\in ]0,r[$,
$\|f\|_{r'}\leq \|f\|_{r}$.
 As a consequence of the assertion (i), we have $$\|f\|_{H^2_{r'}}\lesssim \|f\|_{r'}\leq \|f\|_{r}\lesssim \|f\|_{H^2_{r}},$$
which implies the continuous injection $H_{r}^2 \hookrightarrow H_{r'}^2$.\qed

\smallskip

\begin{lemma}\label{lem_norm-1}
For $r\in \CI$ and $r'\in [0,r[$, and $v\in H_{r}^2$,
we have, for $\tilde r:=\frac{r+r'}{2}$,
$$|v|_{1,r'}\lesssim  \frac{|v|_{0,\tilde r}}{r-r'} \lesssim  \frac{\|v\|_{r}}{(r-r')^{3n+1}}.$$
\end{lemma}

\proof Recall that there is a coordinate chart $\left\{\left(W_i,z^{(i)}\right)\right\}$ of $M_{r_*}$, such that $z^{(i)}\in \Delta_1^n$ and that in one of these charts, $\tilde v$ denotes the expression of $v$.
Recalling the definition \re{rho} of $\rho$ and its properties from \rt{GS-thm}, let us define
\begin{eqnarray*}
\|D(\rho\circ (z^{(i)})^{-1})\|_{0} &:=&	\sup_{\mathfrak z\in z^{(i)}(M_{r_*}\cap W_i)}\sup_{\substack{\zeta\in \mathbb{C}^n,\\|\zeta|\leq 1}}|D(\rho\circ (z^{(i)})^{-1})(\mathfrak z)\zeta|,\\
\delta&:=&  \frac{r-r'}{2\|D(\rho\circ (z^{(i)})^{-1})\|_{0}} \, , \\
z^{(i)}(M_{r'}\cap W_i)_\delta&:=& \left\{z\in\C^n:|z-z^{(i)}(q)|<\delta \ {\rm for \; some} \; q\in M_{r'}\cap W_i\right\}.
\end{eqnarray*}
Assume that $\delta$ is small enough so that the $\delta$-neighborhood of $W_i\cap M_{r'}$ is still in $W_i\cap M_{r_*}$~:
$$\left(z^{(i)}\right)^{-1}\left(z^{(i)}(M_{r'}\cap W_i)_\delta\right)\subset M_{r_*}\cap W_i.$$
Let us devise a Cauchy-like estimate relative to K\"ahler norm. First of all, we can assume that, on the trivialization,
\begin{equation}\label{norm-triv}
|\zeta|_{z,\kappa}^2\simeq \sum_{1\leq i\leq n}|\zeta_i|^2,\qquad \forall \ z\in\Delta_r^n, \quad  \zeta\in \mathbb{C}^n.
\end{equation}
Let $\tilde v$ be a holomorphic vector field on $\Delta_r^n$.
For $z\in \Delta_{r-\delta}^n$ and $\zeta$ in the unit ball of $\mathbb{C}^n$ (i.e., $\sum_{i=1}^n |\zeta_i|^2=1$), we have, through Cauchy-Schwarz inequality, that
\begin{equation}\label{D-tildev}
|D\tilde v(z)\zeta|_{z,\kappa}^2\lesssim \sum_{1\leq i\leq n} \left|\sum_{1\leq k\leq n}\frac{\partial \tilde v_i}{\partial z_k}(z)\zeta_k\right|^2
\lesssim\sum_{1\leq i\leq n}\sum_{1\leq k\leq n} \left|\frac{\partial \tilde v_i}{\partial z_k}(z)\right|^2.
\end{equation}
Since, through Cauchy estimate, we have that
$$\sum_{1\leq k\leq n}\left|\frac{\partial \tilde v_i}{\partial z_k}(z)\right|^2 \leq  \frac{n}{\delta^2} \sup_{w\in \Delta_r^n}| \tilde v_i(w)|^2\leq  \frac{n}{\delta^2} \sup_{w\in \Delta_r^n} \sum_{1\leq i\leq n} | \tilde v_i(w)|^2.$$
Hence, by \re{norm-triv} and (\ref{D-tildev}) we have~:
$$|D\tilde v(z)\zeta|_{z,\kappa}^2\lesssim \frac{n^2}{\delta^2}\sup_{w\in \Delta_r^n}\sum_{1\leq i\leq n} | \tilde v_i(w)|^2\leq \frac{n}{\delta^2}\sup_{w\in \Delta_r^n}|\tilde v|_{w,\kappa}^2.$$
 According to the definition of norm (\ref{norm01r}), we have, by Cauchy estimate, that
\beq\label{maj1}
|v|_{1,r'}=\max_{i} \sup_{q'\in W_i\cap M_{r'}}\|D\tilde v(z^{(i)}(q'))\|_{\kappa}
\lesssim \max_{i} \sup_{q'\in W_i\cap M_{r'}}\sup_{|z-z^{(i)}(q')|=\delta}\frac{|\tilde v(z)|_{\kappa}}{ \delta}.
\eeq
We recall that $q\in M_{r'}$ if and only if $0\leq \rho(q)<r'$.
With $\tilde r=\frac{r+r'}{2}$, let us show that
\beq\label{maj2}
\max_{i} \sup_{q'\in W_i\cap M_{r'}}\sup_{|z-z^{(i)}(q')|=\delta}|\tilde v(z)|_{\kappa} \leq |v|_{0,\tilde r}.
\eeq
Applying Taylor formula of order $1$, we obtain that
\begin{eqnarray*}
\left|\rho((z^{(i)})^{-1}(z))-\rho(q')\right|&\leq&\|D(\rho\circ (z^{(i)})^{-1})\|_{0} \cdot \left|z-z^{(i)}(q')\right|\\
&\leq & \delta  \|D(\rho\circ (z^{(i)})^{-1})\|_0 \ =  \ \frac{r-r'}{2} \ =  \ \tilde r -r'.
\end{eqnarray*}
Hence,
\beq\label{maj3}
\rho((z^{(i)})^{-1}(z))\leq \rho(q')+ \left|\rho((z^{(i)})^{-1}(z))-\rho(q')\right|  \leq r'+ (\tilde r-r')=\tilde r.
\eeq
The first estimate is shown. The second follows from the latter together with \re{supnormequiv-1}.\qed

\subsection{Action on $TM$}\label{actionTM}

By an {\it isometry} of $M$, we mean an analytic diffeomorphism of $M$ which preserves the distance induced by the Riemannian metric $g$.
We recall that the group of isometries $\text{Isom}(M)$ of $M$ is compact \cite{Kob72}[Chap. 2, Theorem 1.2].
As a diffeomorphism, the (element of the) group $\text{Isom}(M)$ acts on the space of sections $L^2(M,TM)$ by push-forward: if $f\in \text{Isom}(M)$ and $v\in L^2(M,TM)$ then $f_*v(x)=Df(f^{-1}(x))v(f^{-1}(x))$ a.e. . This action is {\bf unitary}~:
\beq\label{scal-tg}
\langle f_*v,f_*w\rangle=\int_M \langle f_*v,f_*w\rangle_{g,x}d{\rm vol}(x)=\int_M \langle v,w\rangle_{f^{-1}_*g,f^{-1}(x)} d{\rm vol}(f^{-1}(x))=\langle v,w\rangle .
\eeq
This action of the group of isometry commutes with the Laplacian on $TM$, that is
\beq \label{iso-laplacian}
\Delta_{TM}(f_*v)=f_*\De_{TM}.
\eeq
This is due to the general fact that $f_*\De_{TM,g}(v)=\De_{TM,f_*g}(f_*v)$.
According to Peter-Weyl theorem \cite{Zim90}[p.61], the Hilbert space $L^2(M,TM)$ can be decomposed into an orthogonal sum of finite dimensional subspaces $V_i$ which are irreductible with respect to the action of $\text{Isom}(M)$~:
\beq\label{peter-weyl}
L^2(M,TM)=\bigoplus_{i\geq 0} V_i.
\eeq
In particular, they are invariant $f_*V_i\subset V_i$, for all $f\in \text{Isom}(M)$ and index $i$. Hence, after reordering the indices, each $V_i$ is contained in the eigenspace $E_{\lam_i}$ associated to an eigenvalue $\lam_i$ of the Laplacian $\De_{TM}$.

\begin{lemma}\label{lem_exp-isometry}
	Let $\pi$ be an analytic isometry of the manifold $M$. For $w\in H^2_{r}$, $0<r<r_*$, with $\|w\|_{r}$ sufficiently small, we have
\begin{eqnarray*}
\pi\circ \Exp\{w\} \circ \pi^{-1}=\Exp \left\{D\pi\cdot w\circ\pi^{-1}\right\}= \Exp \left\{\pi_* w\right\}.
\end{eqnarray*}
\end{lemma}
\proof The diffeomorphism $\pi$ extends to a holomorphic isometry of $M_{r}$ w.r.t. K\"ahler metric $\ka$.
 According to \rrem{complex-riemann}, the complex exponential map also defines the geodesics w.r.t. $\kappa$.
 Hence, we have
 $$\pi\circ \Exp_q\{w\}=\Exp_{\pi(q)}\left\{D\pi(q)\cdot w(q)\right\}.\qed$$

%
%


\section{Local rigidity of group action by isometries}

In this section, with some necessary detailed definitions, we will give the precise statement of the main result.

\subsection{Group action and cohomology}\label{subsec_GroupAction}

For a {\it finitely presented group} $G$, let us fix its presentation, i.e., a finite collection $S=\{\gamma_1, \cdots, \gamma_k\}$ of generators and their inverses,
as well as a finite collection of relations $R=\{\CW_1, \cdots, \CW_p\}$, where each $\CW_i$ is a finite word of elements in $S$. In other words, we can view each $\CW_j$ as a word in an alphabet  of the $k$ letters $\{\gamma_l\}_{1\leq l\leq k}$,
$$\CW_j =\gamma_{l^{(j)}_1}\cdots \gamma_{l^{(j)}_{m_j}},\quad 1\leq l^{(j)}_1, \cdots , l^{(j)}_{m_j}\leq k, \quad j=1,\cdots,p.$$
Furthermore, for each $1\leq j\leq p$, we have $\CW_j=e$, the identity of $G$.

For the real analytic compact Riemannian manifold $(M,g)$, let
$$\pi:G\to {\rm Isom}(M,g)\subset {\rm Diff}^\omega(M)$$
be a morphism group, which defines a {\it $G-$action by analytic isometries}.
Let us consider the {\it $G-$action by analytic diffeomorphisms} $\pi_0: G\to {\rm Diff}^\omega(M)$,
 which is considered as a small perturbation of $\pi$. For $\gamma\in S$, we shall write
 \begin{equation}\label{intro-pi}
 \pi_0(\gamma):=\Exp\{P_0(\gamma)\}\circ \pi(\gamma) \;\  {\rm with } \;\  P_0(\gamma)\in \Gamma^\omega=\Gamma^\omega(M,TM),
 \end{equation}
where $\Exp$ denotes the exponential associate to the Riemannian connection. 
As it extends to a holomorphic section on $M_{r_*}$, we consider that $P_0(\gamma)\in H^2_{r_0}$ for some $r_0\in ]0,r_*[$ for all $\gamma\in S$,  which is assumed to be sufficiently small such that the exponential $\Exp\{P_0(\gamma)\}$ is well-defined.

We shall write $C^i(G,\Gamma^\omega)$ for $i-$cochains on $G$ with values in $\Gamma^\omega$.
One can identify $0-$cochains $C^0(G,\Gamma^\omega)$ with $\Gamma^\omega$, $1-$cochains $C^1(G,\Gamma^\omega)$ with maps from $S$ to $\Gamma^\omega$,
 or equivalently $(\Gamma^\omega)^k$, and $2-$cochains $C^2(G,\Gamma^\omega)$ with maps from $R$ to $\Gamma^\omega$, or equivalently $(\Gamma^\omega)^p$.
The $G-$action $\pi$ induces a representation on $\Gamma^\omega$ defined to be
$$\pi(\gamma)_*v :=
\left(D\pi(\gamma)\cdot v\right)\circ \pi\left(\gamma^{-1}\right) ,\quad v\in\Gamma^\omega,$$
which means the differential of $\pi(\gamma)$ evaluated at $\pi\left(\gamma^{-1}\right)$ and applied on $v\left(\pi\left(\gamma^{-1}\right)\right)$.
To this representation, we can associate the Hoshchild complex:
$$\Gamma^\omega \stackrel{d_0}{\longrightarrow} C^1(G,\Gamma^\omega) \stackrel{d_1}{\longrightarrow} C^2(G,\Gamma^\omega)  \stackrel{d_2}{\longrightarrow}  \cdots ,$$
or equivalently,
\begin{equation}\label{HoshchildComplex}
\Gamma^\omega \stackrel{d_0}{\longrightarrow} (\Gamma^\omega)^k  \stackrel{d_1}{\longrightarrow} (\Gamma^\omega)^p  \stackrel{d_2}{\longrightarrow}  \cdots ,
\end{equation}
where the differentials $d_0$ and $d_1$ can be written explicitly as
\begin{itemize}
\item for $v\in\Gamma^\omega$,
\begin{equation}\label{defi_d0}
d_0v=\left(v-\pi(\gamma_l)_*v\right)_{1\leq l \leq  k},
\end{equation}
\item for $v\in (\Gamma^\omega)^k$ and the finite words $\CW_j=\gamma_{l^{(j)}_1}\cdots \gamma_{l^{(j)}_{m_j}}$, $1\leq j\leq p$,
\begin{equation}\label{defi_d1}
	d_1 v=\left(\sum_{1\leq z \leq m_j}\left(D\pi\left(\prod_{i=1}^{z-1}\gamma_{l^{(j)}_i}\right)\cdot v_{l^{(j)}_z}\right)\circ \left(\prod_{i=z}^{m_j}\pi\left(\gamma_{l^{(j)}_i}\right) \right) \right)_{1\leq j\leq p},
\end{equation}
where $\left(D\pi(\gamma)\cdot v\right)\circ \pi(\gamma')$ means the differential of $\pi(\gamma)$ evaluated at $\pi(\gamma')$ and applied on $v(\pi(\gamma'))$, and for $z=1$, we interpret the empty product ``$\prod_{i=1}^{z-1}$" as the identity $e\in G$.
\end{itemize}
It is obvious that, for $r > 0$, the linear operators $d_0$ and $d_1$ are bounded on $\Gamma^\om_r$ and $(\Gamma^\om_r)^k$ respectively.

Following the idea of Lombardi-Stolovitch \cite{ LS10} developed in the context of germs of holomorphic vector fields at a singular point, let us consider the {\it self-adjoint box} operator
\begin{eqnarray}
\square:  C^1(G,\Gamma^\omega) &\to& C^1(G,\Gamma^\omega)\label{def-box}\\
v &\mapsto& \square v:=(d_0\circ d_0^*+d_1^*\circ  d_1)v  \   ,\nonumber
\end{eqnarray}
where $d_0^*$ and $d_1^*$ denote the adjoints of $d_0$ and $d_1$ respectively, with respect to the scalar product from $L^2(M,TM)$ induced on $(L^2(M,TM))^k$ and  $(L^2(M,TM))^p$ respectively.

By definition, both $d_0$ and $d_1$ are defined by the $G-$action by isometries on $L^2$ vector fields.
According to Peter-Weyl theorem mentioned above, the decomposition into irreducible finite-dimensional subspace (\ref{peter-weyl}) gives rise to complexes:
\begin{equation}\label{chain-invariant-d0d1}
V_i\stackrel{d_0}{\longrightarrow}  V_i^k \stackrel{d_1}{\longrightarrow}   V_i^p,\qquad  \forall  \  i\in\N.
\end{equation}
Recalling the decomposition (\ref{decomp_lambdas}), let $\PP_i$ be the projection onto $V_i$ or $V_i^k$ (depending on the context), i.e.,
\begin{equation}\label{Proj_i}
\PP_i u=\left\{\begin{array}{cl}
  \displaystyle \sum_{j\in I_i} u_j {\bf e}_j, &   \displaystyle  u= \sum_{j\in \N} u_j {\bf e}_j\in L^2(M,TM)\\[2mm]
  \displaystyle \left(\sum_{j\in I_i} u_{l,j} {\bf e}_j\right)_{1\leq l\leq k}, &  \displaystyle  u= \left(\sum_{j\in \N} u_{l,j}  {\bf e}_j\right)_{1\leq l\leq k}\in (L^2(M,TM))^k
\end{array}\right.  .\end{equation}
Then the box operator $\square_i:=\square\circ \PP_i$ is a self-adjoint operator of finite-dimensional vector space.
Hence it is diagonalizable and its non-zero eigenvalues $\mu_{i,1},\cdots,\mu_{i,K_i}$ (if existing) are all positive.
In that case, let us define $\mu_i$ to be the smallest non-zero eigenvalue of $\square_i$, i.e.,
\begin{equation}\label{mu_i}
\mu_i:=\min_{1\leq j \leq K_i}\mu_{i,j}.
\end{equation}
Through elementary properties of Hilbert space, we have
\begin{equation}\label{decomVi_ImKer}
V_i^k  = {\rm Im}\square_i \bigoplus {\rm Ker}\square_i.
\end{equation}
In view of Section \ref{actionTM}, $V_i$ is contained in the eigenspace $E_{\lambda_i}$ of $|\Delta_{TM}|^{\frac12}$ associated with the eigenvalue $\lambda_i$.
Since $\square_i$ is diagonalizable on $V_i^k$, let $\{{\CE}_{i,1}, \cdots, {\CE}_{i,K_i}\}\subset E_{\lambda_i}^k$ be an orthonormal basis of its eigenvectors associated with the non-zero eigenvalues $\{\mu_{i,1}, \cdots, \mu_{i,K_i}\}$ of $\square_i$, and let $\{\mathcal{N}_{i,1},\ldots ,\mathcal{N}_{i,J_i}\}$ be an orthonormal basis of ${\rm Ker}\square_i$.
Hence every $u\in V_i^k$ can be decomposed along these eigenvectors, i.e.,
\begin{equation}\label{u_EquiForms}
\left(\sum_{j\in I_i} u_{l,j} {\bf e}_j\right)_{1\leq l\leq k} = u=\sum_{1\leq j\leq K_i}\tilde u_{i,j} \CE_{i,j}+ \sum_{1\leq j\leq J_i}\breve u_{i,j} \CN_{i,j} ,
\end{equation}
which implies that
\begin{equation}\label{u_EquiForms-norm}
\left(\sum_{1\leq l\leq k}\sum_{j\in I_i} \left|u_{l,j} \right|^2 \right)^{\frac12} = \left\|u\right\|_{L^2}= \left(\sum_{1\leq j\leq K_i}\left|\tilde u_{i,j}\right|^2+\sum_{1\leq j\leq J_i}\left|\breve u_{i,j}\right|^2\right)^{\frac12}.
\end{equation}

\begin{lemma}\label{lem-H1}
Given the $G-$action $\pi$ on $M$ as above, we have
$$H^1(G, V_i):=\Ker \,(d_{1}\circ\PP_i )/ \IM \, (d_{0}\circ\PP_i )\approxeq {\rm Ker}\square_i.$$
\end{lemma}
\proof Given $f\in V^k_i$, we have that
	$$
	\la \square_i f, f\ra = \|d_0^*f\|^2+\|d_1f\|^2
	$$
	Hence, $\square_i f=0$ iff $d_0^*f=0$ and $d_1f=0$. The latter means that $f$ is a 1-cocycle, i.e. $f\in Z^1(G,V_i)$ . The former means that $f$ is orthogonal to $\text{Im}\,d_0$. Hence, $f$ belongs to a space isomorphic to $Z^1(G,V_i)/\text{Im}(d_{0}\circ\PP_i )=H^1(G, V_i)$. 
\qed

\begin{lemma}\label{actionsurH}
For $0<r<r_*$, the $G-$action by isometries $\pi$ acts on $H_r^2$, and the Hochschild complex \re{HoshchildComplex} gives rise to a complex of Hilbert spaces~:
	\begin{equation}\label{HoshchildComplex-hardy}
		H_r^2 \stackrel{ d_0}{\longrightarrow} (H_r^2)^k  \stackrel{d_1}{\longrightarrow} (H_r^2)^p  \stackrel{ d_2}{\longrightarrow}  \cdots .
	\end{equation}
Moreover, the linear operators $d_0$ and $d_1$ are uniformly bounded on $H^2_r$ and $(H^2_r)^k$ respectively, w.r.t. $0<r<r_*$.
\end{lemma}
\proof
	For $1\leq l\leq k$, $\pi(\gamma_l)$ is an analytic isometry of $M$. Since the K\"alher metric $\kappa$ extends to $M_{r}$ the Riemannian metric $g$ of $M$, $\pi(\gamma_l)$ extends to a holomorphic isometry of $M_{r}$.  Let us show that $\pi(\gamma_l)$ leaves invariant each boundary $\partial M_{r'}$, $0<r'<r$. Indeed, according to \cite{klimek}[4. p.16] and \re{MA},
		$$\det\left(\frac{\partial^2 \sqrt{\rho}\circ \pi(\gamma_l)}{\partial z_i\partial\bar z_j}\right)=|\det D\pi|^2\det\left(\frac{\partial^2 \sqrt{\rho}}{\partial z_i\partial\bar z_j}\right)\circ \pi(\gamma_l)=0.$$
		On the other hand, according to \cite{klimek}[(1.4.4)], we have
		$$
		\la\mathcal{L}({ \rho}\circ \pi(\gamma_l))(a)b,c\ra= \la\mathcal{L}({ \rho})(\pi(\gamma_l)(a))D\pi(\gamma_l)(a)b,D\pi(\gamma_l)(a)c\ra,
		$$
		where $a\in M_r$, $ b,c\in T^{(1,0)}M_r$ and
		$$
		\la\mathcal{L}(u(a)b,c\ra:=\sum_{1\leq i,j\leq n}\frac{\partial^2 u}{\partial z_i\partial \bar z_j}(a)b_i\bar c_j,\quad b:=\sum_{1\leq i\leq n} b_i\frac{\partial}{\partial z_i}, \quad c:=\sum_{1\leq i\leq n} c_i\frac{\partial}{\partial z_i} .
		$$
Since $\pi(\gamma_l)$ is an isometry of both $M$ and $M_{r'}$, recalling the formula of the K\"ahler metric \re{kappa} defined through $\rho$, the previous equality shows that we also have an isometric embedding of $M$ into $M_{r'}$ with the K\"ahler metric associated to the real analytic potential $\rho\circ\pi(\gamma_l)$.
By uniqueness of the solution from \rt{GS-thm}, we have $\rho\circ \pi(\gamma_l)= \rho$. Therefore, $\pi(\gamma_l)$ leaves invariant the level sets of $\rho$. It is also an isometry of $\partial M_{r'}$. As a consequence, according to \re{h2product}, we have for $f, h \in  H_{r'}^2$,
\begin{eqnarray*}
\la \pi(\gamma_l)_*f,\pi(\gamma_l)_*h \ra_{H_{r'}^2}&=&\int_{\partial M_{r'}}\la \pi(\gamma_l)_*f(q),\pi(\gamma_l)_*h(q)\ra_{q,\ka}d\mu_{r}(q)\\
&=&\int_{\partial M_{r'}}\la f,h\ra_{\pi^{-1}(q),\ka}d\mu_{r'}(q)=\int_{\partial M_{r'}}\la f,h\ra_{q,\ka}d\mu_{r'}(\pi^{-1}(q)).
\end{eqnarray*}
Hence, $\pi(\gamma_l)$ acts on the Hardy space $H_r^2$ and
 $\| \pi(\gamma_l)_*f\|_{H_{r'}^2}\lesssim \|f\|_{H_{r'}^2}$.

As the expressions \re{defi_d0} and \re{defi_d1} involve only terms of the form $(D\pi(\gamma)v)$ evaluated at some $\pi(\tilde \gamma)$, with $v$ a holomorphic vector field, $\left. d_0\right|_{H^2_r}$ and $\left. d_1\right|_{(H^2_r)^k}$ take their values in $(H^2_r)^k$ and $(H^2_r)^p$ respectively.\qed

 \smallskip

From now on, let $\|d_0\|_{r}$ and $\|d_1\|_{r}$ denote the corresponding operator norms on $H^2_r$ and $(H^2_r)^k$ respectively.

\subsection{Precise statement of main theorems}

 The Diophantine property of the self-adjoint box operator $\square$ (hence of the $G-$action by isometries $\pi$) mentioned in Theorem \ref{thmmain} and \ref{thm-geom0} is precisely defined through the asymptotic lower bound of the non-zero eigenvalues of $\square_i$ w.r.t. $i\in\N$. Recall $\mu_i$ introduced in (\ref{mu_i}), we define the Diophantine condition for $\pi$ as follows.

\begin{defi}\label{defi-dio} The $G-$action by isometries $\pi$ is called {\bf Diophantine} if there exists non-negative numbers $\sigma$, $\tau$ such that the operator $\square$ satisfies
\begin{equation}\label{diophantine}
\mu_i\geq \frac{\sigma}{(1+\lambda_i)^\tau},\qquad \forall \ i\in \N.
\end{equation}
\end{defi}


Recall the $G-$action by analytic diffeomorphisms $\pi=\Exp\{P_0\}\circ \pi$ given in (\ref{intro-pi}),
which is considered as a small perturbation of the $G-$action by analytic isometries $\pi$, with $P_0:S\to H^2_{r_0}$ sufficiently small in the sense that
\begin{equation}\label{norm-Gr}
\|P_0\|_{S,r_0}:=\left(\sum_{1\leq l \leq k}\|P_0(\gamma_l)\|^2_{r_0}\right)^{\frac12}
\end{equation}
is sufficiently small.
As a more general notation, for $u=u(\gamma)$ defined on $S$ equipped with some norm $\|\cdot\|_*$ (or sometimes written as $|\cdot|_*$), we define
\begin{equation}\label{norm-G}
\|u\|_{S,*}=\left(\sum_{1\leq l \leq k}\|u(\gamma_l)\|^2_*\right)^{\frac12}.
\end{equation}
Recalling the definition of norms \re{norm-moser}, we define~:

\begin{defi}\label{defi_formal_conj}
Let $\pi'$ be a $G-$action by analytic diffeomorphisms on $M$. Given $0<\zeta<1$, $\pi'$ is said to be {\bf $\zeta-$formally conjugate} to $\pi$ on $M$, and we write $\pi'\in \CF^M_{\pi,\zeta}$, if, for any $\varepsilon>0$, there exists $y^\varepsilon \in \Ga^\om$ with $\|y^\varepsilon\|_{1,M}<\zeta$ such that, for every $\gamma\in S$,
\begin{equation}\label{eq_form_conj}
\Exp\{y^\varepsilon\}^{-1} \circ \pi'(\gamma)\circ  \Exp\{y^\varepsilon\}=\Exp\{z^\varepsilon(\gamma)\} \circ \pi(\gamma),
\end{equation}
for some $z^\varepsilon:S\to \Ga^\om$ satisfying $\|z^\varepsilon\|_{S,1,M}<\varepsilon$.

\end{defi}

\begin{remark} The definition of $\pi'\in \CF^M_{\pi,\zeta}$ means that $\pi'$ can be conjugated to $\pi$, as close as one wishes, by means of near-identity analytic diffeomorphisms.
Note that the nature of $\{y^\varepsilon\}\subset\Ga^\om$ (i.e., convergence w.r.t. some topology as $\varepsilon\to 0$) is not mentioned in the above definition, which makes the essential difference from the analytic conjugation of $\pi'$ to $\pi$.
\end{remark}


The precise statements of Theorem \ref{thmmain} and \ref{thm-geom0} are given as follows.

\begin{thm}\label{thm_conj} 
Let $\pi$ be a Diophantine $G-$action by analytic isometries on the analytic Riemannian manifold $M$. Assume that 
$\pi$ satisfies $\dim\Ker \square<\infty$.

Let $\pi_0\in \CF^{M}_{\pi,\zeta}$ be a $G-$action by analytic diffeomorphisms on $M$, and 
$$\pi_0(\gamma)=\Exp\{P_0(\gamma)\}\circ \pi(\gamma),\quad \gamma\in S,$$ with $P_0:S\to H_{r_0}^2$ for some $r_0\in ]0, r_*[$. 
If $\|P_0\|_{S,r_0}:=\varepsilon_0$ sufficiently small, and $0<\zeta<|\ln(\varepsilon_0)|^{-\frac13}$, then $\pi_0\in \CF^{M}_{\pi,\zeta}$ is analytically conjugate to $\pi$.
\end{thm}

\begin{thm} \label{thm-geom}
Let $\pi$ be a Diophantine $G-$action by analytic isometries on $M$ as above. Assume that $H^1(G,V_i)=0$ for all $i\in \N$. Then any small enough analytic perturbation $\pi_0$ of $\pi$ is analytically conjugate to $\pi$.
\end{thm}

\subsection{Examples of local analytic rigidity results}\label{exemple-sect} 


 Before stating the corollaries of Theorem \ref{thm_conj} and \ref{thm-geom}, let us introduce another definition of Diophantine property of group actions, given by Dolgopyat \cite{dolgopyat} who defined a ``small divisors condition" for the subset of the group $G$.

\begin{defi}(Dolgopyat \cite{dolgopyat}[Appendix A])\label{dioph-dolg}
Given a finitely presented group $G$ acting transitively on a compact manifold $M$ by isometries $\pi$, the subset $\mathcal S\subset G$ is said to be a {\bf Diophantine subset for $\pi$} if there exist constants $C>0$ and $\tau\geq 0$ such that,  for all $i\in \N$ with $\lambda_i\neq 0$ and all $f_i\in V_i$, there exists $\gamma\in {\mathcal{S}}$ (depending on $f_i$) such that 
	\begin{equation}\label{dioph-cond}
\|f_i -\pi(\gamma )_*f_i\|_{L^2} \geq\frac{C}{\lambda_i^{\tau}}\|f_i\|_{L^2}.
	\end{equation}
\end{defi}
 We recall that $R=\{\CW_1,\cdots, \CW_p\}$ denotes the set of relators of $G$. Each $\CW_j= \gamma_{l^{(j)}_1}\cdots\gamma_{l_{m_j}^{(j)}}$ is a word in the alphabet $S=\{\gamma_{1},\cdots,\gamma_{k}\}$ composed of the generators and their inverse of $G$. Since $\pi(\CW_j)={\rm Id}$, according to \re{defi_d1}, we have
\begin{eqnarray}
	d_1 v&=&\left(\sum_{1\leq q \leq m_j}\left(D\pi\left(\prod_{i=1}^{q-1}\gamma_{l^{(j)}_i}\right)\cdot v_{l^{(j)}_q}\right)\circ \left(\pi\left(\prod_{i=1}^{q-1}\gamma_{l^{(j)}_i}\right)^{-1} \right) \right)_{1\leq j\leq p},\nonumber\\
	&=&\left(\sum_{1\leq q \leq m_j}\pi\left(\prod_{i=1}^{q-1}\gamma_{l^{(j)}_i}\right)_* v_{l^{(j)}_q} \right)_{1\leq j\leq p},\qquad v\in (L^2(M,TM))^k.\label{defi_d1b}
\end{eqnarray}
 As the expression $d_1$ is defined by the Fox derivatives of the relators \cite{birman, curran}, we  use the convention that $\gamma_{i+\nu}=\gamma^{-1}_i$ and $v_{i+\nu}:=-\pi(\gamma_i^{-1})_*v_i$, $1\leq i\leq \nu$, whenever $\{\gamma_1,\cdots, \gamma_\nu\}$ denotes a set of generators of $G$.
	
 Let 
$
\CS_j:=\left\{\gamma_{l^{(j)}_1}\cdots \gamma_{l^{(j)}_q}: 1\leq q\leq m_j\right\}
$
be the set of subwords of $\CW_j$, and let
$$
\CS_{j,i}:=\left\{w'\in \CS_j\cup\{e\}: w'\gamma_{i}\in \CS_j\right\},\qquad i=1,\cdots, k,
$$
be the set of words $w'$ of $\CS_j\cup\{e\}$ such that $w'\gamma_{i}$ is still a word of $\CS_j$.
For $(d_1 v)_j$, the $j$-th coordinate of $d_1v$, we have
\begin{equation}\label{d1-new}
(d_1 v)_j=\sum_{1\leq q \leq m_j}\pi\left(\prod_{i=1}^{q-1}\gamma_{l^{(j)}_i}\right)_* v_{l^{(j)}_q}= \sum_{l=1}^k \sum_{w\in \CS_{j,l}}\pi(w)_*v_l
\end{equation}
with the convention that the corresponding sum equal $0$ if $\CS_{j,l}=\emptyset$. 
For $W=(W_1,\cdots, W_p)\in (L^2(M,TM))^p$, using the scalar product on $(L^2(M,TM))^p$, we have
\begin{eqnarray*}
	\la d_1v, W\ra &=& \sum_{j=1}^p\la (d_1 v)_j, W_j\ra\\
	&=& \sum_{j=1}^p \left\la \sum_{l=1}^k\sum_{w\in \CS_{j,l}}\pi(w)_*v_l, W_j\right\ra
	\ = \  \sum_{l=1}^k \left\la v_l, \sum_{j=1}^p\sum_{w\in \CS_{j,l}}\pi(w^{-1})_*W_j\right\ra,
\end{eqnarray*}
which implies that
\beq\label{d1*}
d_1^*W:= \left( \sum_{j=1}^p \sum_{w\in \CS_{j,l}}\pi(w^{-1})_*W_j \right)_{1\leq l \leq k}.
\eeq
A substitution of \re{d1-new} into \re{d1*} shows that 
$$
(d_1^*d_1)v= \left(\sum_{m=1}^k\sum_{j=1}^p\sum_{w\in\CS_{j,l}}\sum_{\tilde w\in \CS_{j,m}}\pi(w^{-1}\tilde w)_*v_m\right)_{1\leq l\leq k} .
$$
In this sense, $d_1^*d_1$ can be regarded as a self-adjoint $k\times k$-matrix whose $l$-row and $m$-column term is
$$
(d_1^*d_1)_{l,m}:=\sum_{j=1}^p \sum_{w\in \CS_{j,l}}\sum_{\tilde w\in \CS_{j,m}}\pi(w^{-1}\tilde w)_* .
$$
It is easy to verify that, for any $i\in \N^*$, all eigenvalues of $\left.(d_1^*d_1)\right|_{V_i^k}$ are non-negative.


\begin{defi}\label{GDpi} The finitely presented group $G$ is called {\bf Diophantine w.r.t $\pi$} if 
\begin{enumerate}
\item\label{DCG} the set of generators of $G$ is a Diophantine subset for $\pi$.
\item\label{DCR} the set of relators of $G$ is Diophantine in the sense that there exist $C',\tau'>0$ such that, for all $i\in \N^*$, the smallest non-zero eigenvalue of $\left.(d_1^*d_1)\right|_{V_i^k}$ is bounded from below by $C'\lambda_i^{-\tau'}$.
\end{enumerate}
\end{defi}

\begin{prop}\label{DG-D} Let $\pi$ be an $G-$action by isometries as above. If the  group $G$ is Diophantine w.r.t. $\pi$ 
	, then $\pi$ is Diophantine in the sense of \rd{defi-dio}.
\end{prop}
\proof Let $\{\gamma_1,\cdots,\gamma_\nu\}$ be the set of generators of $G$. By setting $\gamma_{i+\nu}:=\gamma_{i}^{-1}$, we have $S=\{\gamma_j\}_{1\leq j\leq k}$ with $k=2\nu$. In view of \cite{dolgopyat}[Proposition A.1-(b)], $\{\gamma_1,\cdots,\gamma_\nu\}$ is a Diophantine subset for the $G-$action $\pi$ if and only if $S$ is a Diophantine subset for $\pi$.

For $1\leq j\leq k$,  let us define $L_i$ on $L^2(M,TM)$ by 
$$L_ju :=u -\pi(\gamma_j)_*u, \qquad \forall \ u\in V_i, \quad i\in \N.$$
Then, for $u\in V_i$, we have
 $$\|d_0u\|_{L^2}^2=\sum_{j=1}^k\|L_ju\|_{L^2}^2.$$ 
 If $\{\gamma_1,\cdots,\gamma_\nu\}$ is a Diophantine subset for $\pi$ (so is $S$), then for any non-zero $u\in V_i$ with $\lambda_i > 0$, there exists $j$
such that $\|L_ju\|_{L^2}\geq \frac{C}{\lambda_i^{\tau}}\|u\|_{L^2}$, and hence 
$$\|d_0u\|_{L^2}\geq \frac{C}{\lambda_i^{\tau}}\|u\|_{L^2}.$$
In particular, for all $v\in V_i^k$ such that $d_0^*v\neq 0$, we have 
$$\|(d_0\circ d_0^*)v\|_{L^2}\geq \frac{C}{\lambda_i^{\tau}}\|d_0^*v\|_{L^2}.$$ 
Let $v$ be an eigenvector of $\square_i$ associated to the positive eigenvalue $\mu$~: 
\begin{equation}\label{ev-box}
\square v =(d_0\circ d_0^*)v+(d_1^*\circ d_1)v=\mu v.
\end{equation}
As the operators $d_0\circ d_0^*$ and $d_1^*\circ d_1$ take values in mutual orthogonal spaces ${\rm Ker}d_1$ and ${\rm Im}d_1^*$ respectively, and $d_1\circ d_0=0$, the equality (\ref{ev-box}) can be decomposed into 
$$(d_0\circ d_0^*)v^{\rm Ker}=\mu v^{\rm Ker},\quad (d_1^*\circ d_1)v^{\rm Im}=\mu v^{\rm Im},$$ 
where $v=v^{\rm Ker}+v^{\rm Im}$ with $v^{\rm Ker}\in{\rm Ker}d_1$, $v^{\rm Im}\in {\rm Im}d_1^*$. If $d_0^*v^{\rm Ker}\neq 0$, then we have $$\mu\|v^{\rm Ker}\|_{L^2}=\|(d_0\circ d_0^*)v^{\rm Ker}\|_{L^2}\geq \frac{C}{\lambda_i^{\tau}}\|d_0^*v^{\rm Ker}\|_{L^2}.$$ On the other hand, since
$$
\|d_0^* v^{\rm Ker}\|_{L^2}^2=\left\langle (d_0\circ d_0^*) v^{\rm Ker}, v^{\rm Ker}\right\rangle=\mu \|v^{\rm Ker}\|_{L^2}^2,
$$
we have
$\sqrt{\mu}\geq \frac{C}{\lambda_i^{\tau}}$. If $d_0^*v^{\rm Ker}= 0$, then $v^{\rm Im}\neq 0$. According to \rd{GDpi}(\ref{DCR}), $\mu\geq \frac{C'}{\lambda_i^{\tau'}}$.
Hence the inequality (\ref{diophantine}) in Definition \ref{defi-dio} holds for all $i\in\N^*$. As for $\lambda_0=0$, since $\square_0$ is finite-dimensional, its non-zero eigenvalues (if existing) are bounded from below by some constant.
\qed


\medskip


As a corollary of Theorem \ref{thm-geom}, we obtain ~:
\begin{cor}
Let $G$ be a discrete group with Kazhdan's property (T) and let $\pi$ be $G-$action by analytic isometries on $M$ as above having a Diophantine set of relators in the sense of \rd{GDpi}\re{DCR}. Then any small enough analytic perturbation $\pi_0$ of $\pi$ 
 is analytically conjugate to $\pi$.
\end{cor}
 \begin{proof}
	According \cite{BHV}[Theorem 3.2.1], for all $i\in \N$, $H^1(G,V_i)=0$. On the other hand, $(T)$ condition implies that 
	there exist constants $C>0$ such that,  for all $i\in \N$, and all $f_i\in V_i$, there exists $\gamma\in S$ (depending on $f_i$) such that $\|f_i -\pi(\gamma )_*f_i\|\geq C\|f_i\|$. Hence, $S$ is Diophantine with $\tau=0$. So we can apply \rp{DG-D} and Theorem \ref{thm-geom}.
	\end{proof}
In the smooth category, this result is due to Fisher-Margulis theorem \cite{margulis-fisher-annals}[Theorrem 1.3] (see also \cite{Fisher}[Theorem 1.2]) and holds without the Diophantine assumption. The following result might be known to experts although we haven't found a reference.
\begin{cor} Let $G$ be a compact connected semi-simple Lie group.
	Assume that the smallest closed subgroup of $G$ containing $S$ is $G$. Let $\pi$ be $G-$action by analytic isometries on $M$ as above having a Diophantine set of relators in the sense of \rd{GDpi}\re{DCR}. Let us assume that $H^1(G,V_i)=0$ for all $i\in \N$. Then any small enough analytic perturbation $\pi_0$ of $\pi$ is analytically conjugate to $\pi$.
\end{cor}
\proof
	According to  \cite{dewitt}[Proposition 7], $S$ is Diophantine in the sense of Definition \ref{dioph-dolg}. So we can apply Theorem \ref{thm-geom}.
\qed

\medskip

The following corollary is the analytic version of \cite{Fisher}[Theorem 1.3]~:
\begin{cor}
Let $\Gamma$ be an irreducible lattice in a semi-simple Lie group with rank at least $2$. Let $\pi$ be an action of $\Gamma$ by analytic isometries on a compact real analytic Riemannian manifold $M$ as above  having a Diophantine set of relators in the sense of \rd{GDpi}\re{DCR}. Then any small enough analytic perturbation $\pi_0$ is analytically conjugate to $\pi$.
\end{cor}
\proof
According to Margulis \cite{margulis-book}[introduction, theorem 3], for all $i\in \N$, $H^1(\Gamma, V_i)=0$. Furthermore, $\overline{ \langle \Gamma\rangle}$ is a semi-simple Lie group. According to Dolgopyat \cite{dolgopyat}[Theorem A.3] and \cite{dewitt}[Proposition 7], $\Gamma$ is Diophantine in the sense of Definition \ref{dioph-dolg}. We can then apply \rt{thm-geom}.
\qed
\medskip


 We don't know in general if the set of relators of $G$ is Diophantine whenever its set of generators is. However, the following result shows that it is so in the abelian case.
\begin{cor}
If $G$ is abelian and has a Diophantine set of generators w.r.t. $\pi$, then the $G-$action $\pi$ is also Diophantine.
\end{cor}

\begin{proof}
Let ${\mathcal G}=\{\gamma_1,\cdots,\gamma_k\}$ be a set a generator of $G$. 
Let us set $\gamma_{i+k}:=\gamma_{i}^{-1}$.
Let $U\in L^2(M, TM)$ and $V=(v_1,\cdots, v_k)\in L^2(M, TM)^k$. We have, by definition,
$$d_0U:=\begin{pmatrix}U-\pi(\gamma_i)_*U\end{pmatrix}_{1\leq i\leq k}=\begin{pmatrix}
L_i(U)
\end{pmatrix}_{1\leq i\leq k},\quad d_0^*V:=\sum_{i=1}^kL_i^*(v_i)=\sum_{i=1}^k(v_i-\pi(\gamma_i^{-1})_*v_i).
$$
Assuming $G$ is abelian, so the relations are defined by the words ${\mathcal W}_{i,j}:=\gamma_i\gamma_j\gamma_{i+k}\gamma_{j+k}$ with $1\leq i<j\leq k$. Expressing $d_1$ in term of Fox derivatives \cite{birman, curran}, we have
$$d_1V= \left(v_i+\pi(\gamma_i)_*v_j+\pi(\gamma_i\gamma_j)_*v_{i+k}+\pi(\gamma_i\gamma_j\gamma_{i+k})_*v_{j+k}\right)_{i,j}.$$
together with $v_{i+k}:=-\pi(\gamma_i^{-1})_*v_i$,
we have
\begin{eqnarray*}
d_1V&=&\left(v_i+\pi(\gamma_i)_*v_j-\pi(\gamma_i\gamma_j\gamma_i^{-1})_*v_{i}-\pi(\gamma_i\gamma_j\gamma_{i}^{-1}\gamma_{j}^{-1})_*v_{j}\right)_{i,j}\\
&=& \left((v_i-\pi(\gamma_j)_*v_i)-(v_j-\pi(\gamma_i)_*v_j)\right)_{i,j}\\
&=& \left(L_j(v_i)-L_i(v_j)\right)_{i,j}.\end{eqnarray*}
Hence, we regard $d_1$ as an operator from $(L^2(M, TM))^k$ to ${\rm so}(k)$.
For $W=(W_{i,j})_{1\leq i,j \leq k}\in {\rm so}(k)$, i.e., $W_{i,j}=-W_{j,i}$, we have
	$$
		\langle d_1V, W \rangle = \sum_{i<j}\langle L_jv_i-L_iv_j,W_{i,j}\rangle
		=  \sum_{i=1}^k \left\la v_i, \sum_{j=1}^k L_j^*W_{i,j}\right\ra.
$$
Since $L_iL_j^*=L_j^*L_i$, we have 
$$
d_1^*d_1V=\left(\CL(v_i)-\sum_{j=1}^k L_iL_j^*v_j\right)_{1\leq i\leq k},\quad \CL(v):=\sum_{j=1}^kL_j^*L_jv.
$$
Then we obtain that 
	\begin{equation}\label{box-abelien}
		\square V= (d_0d_0^*+d_1^*d_1) V = \left(\left(\sum_{j=1}^ kL_j^*L_j\right)v_i\right)_{ 1\leq i\leq k}.
	\end{equation}
	
Let us restrict to $V_i^k$, where $V_i$ is the irreducible finite-dimensional obtained by Peter-Weyl theorem, contained in the $\lambda_i$-eigenspace of $\left|\Delta_{TM}\right|^\frac12$.
Recall that the set of generators ${\mathcal G}$ is Diophantine in the sense of \rd{dioph-dolg}, i.e., for any $U\in V_i$, there exists $1\leq j(U)\leq k$ such that 
\begin{equation}\label{dioph-abelien}
\|U-\pi(\gamma_{j(U)})_*U\|\geq \frac{C}{\lambda_i^{\tau}}\|U\|.
\end{equation}
According to (\ref{box-abelien}), we have 
	$$
	\|\square V\|^2=\sum_{i=1}^k\left\|\left(\sum_{j=1}^k L_j^*L_j\right)v_i\right\|^2,
	$$
and, through Cauchy-Schwartz inequality, 
	$$
	\sum_{j=1}^k\|L_j v_i\|^2=\left\la \left(\sum_{j=1}^kL_j^*L_j\right)v_i, v_i\right\ra \leq \left\|\left(\sum_{j=1}^kL_j^*L_j\right)v_i\right\|\|v_i\|.
	$$
	According to $(\ref{dioph-abelien})$, for each $\ell$, there exists $j$ such that $\|L_jv_{\ell}\|\geq \frac{C}{\lambda_i^{\tau}}\|v_{\ell}\|$.
	Therefore,
	$$
	\left\|\left(\sum_{j=1}^kL_j^*L_j\right)v_{\ell}\right\|\geq \frac{C^2}{\lambda_i^{2\tau}}\|v_{\ell}\|,\qquad \ell=1,\cdots , k, 
	$$
	so that
	$$
	\|\square V\|\geq \frac{C^2}{\lambda_i^{2\tau}}\left(\sum_{\ell=1}^k\|v_{\ell}\|^2\right)^\frac12.
	$$ 
	Hence, the $G-$action $\pi$ is Diophantine in the sense of \rd{defi-dio}.
\end{proof}

Theorem \ref{thm-tore0} can be seen as a corollary of Theorem \ref{thm_conj}.

\noindent{\it Proof of Theorem \ref{thm-tore0}.}
Condition \re{dioph-cond} reads: there exist $c ,\tau>0$, such that for all $({\bf k},l)\in \mathbb{Z}^d\times \mathbb{Z}\setminus\{0\}$,
$$\max_{1\leq i\leq m}|\la {\bf k},\al_i\ra-l|\geq \frac{c}{|{\bf k}|^{\tau}}.$$
As the only way to have $|\la {\bf k},\al_i\ra-l|=0$ for all $i$ is to have $({\bf k},l)=(0,0)\in\mathbb{Z}^{d+1} $, we have $\text{Ker}\,\square_j=0$ if and only if $j\in I_0$ (that is the constant term in the Fourier expansion). Hence $\text{Ker}\,\square$ is finite dimensional. According to Petkovic's theorem\cite{petko-tore}[Theorem 5], the family $\{\pi_0(e_i)\}$ is simultaneously {\it smoothly} conjugate to $\{\pi(e_i)\}$. Hence, $\pi_0$ formally conjugate to $\pi$. 
Therefore, according to \rt{thm_conj}, the family $\{\pi_0(e_i)\}$ is simultaneously {\it analytically} conjugate to $\{\pi(e_i)\}$.
\qed

\section{Properties of box operator $\square$}\label{section-box}

In this section, let us derive some properties related to the box operator $\square$, as well as that of the ingredients $d_0$ and $d_1$. These properties will be applied in the next section for the KAM scheme.

\subsection{$d_0$ and $d_1$} Recall the definitions of linear operators $d_0$ and $d_1$ in (\ref{defi_d0}) and (\ref{defi_d1}).

\begin{lemma}\label{lemma_d0} The operators $d_0^*$ and $d_1^*$ are uniformly bounded on $(H_r^2)^k$ and $(H_r^2)^p$ respectively w.r.t. $r\in[0,r_*[$.
\end{lemma}
\proof For any $U=(U_l)_{1\leq l\leq k}\in  (H_r^2)^k$, we have $\left.U\right|_{M}\in (L^2(M,TM))^k$, and
$$ \la \pi(\gamma_l)_* u, U_l \ra =  \la \pi(\gamma_l)_* u, \pi(\gamma_l)_* (\pi(\gamma_l)_*^{-1} U_l) \ra =  \la u, \pi(\gamma_l)_*^{-1}U_l \ra, \quad u\in L^2(M,TM),$$
where $\la\cdot,\cdot\ra$ denotes the scalar product on $L^2(M,TM)$ or on $(L^2(M,TM))^k$, depending on the context.
As a consequence, we have, for $d_0$ defined in (\ref{defi_d0}),
$$\la d_0 u, U \ra = \sum_{l=1}^k \left\la u-\pi(\gamma_l)_*u, U_l \right\ra= \left\la u, \sum_{l=1}^k \left(U_l-\pi(\gamma_l)_*^{-1} U_l\right) \right\ra. $$
Therefore, we have that
\begin{equation}\label{d_0*v}
d_0^* U=  \sum_{l=1}^k \left(U_l-\pi(\gamma_l)_*^{-1} U_l\right).
\end{equation}
%
According to \rl{actionsurH}, we have
\beq\label{boundd0*}\|d_0^* U\|_r\lesssim \sum_{i=1}^k \|U_i\|_r\lesssim \left(\sum_{i=1}^k \|U_i\|_r^2\right)^{\frac12}.\eeq
Hence, $d_0^*$ is uniformly bounded on $(H_r^2)^k$.
Noting that
$$
\left(D\pi\left(\prod_{i=1}^{z-1}\gamma_{l^{(j)}_i}\right)\cdot v_{l^{(j)}_z}\right)\circ \left(\prod_{i=z}^{m_j}\pi\left(\gamma_{l^{(j)}_i}\right) \right)= \left(\pi\left(\prod_{i=1}^{z-1}\gamma_{l^{(j)}_i}\right)_* v_{l^{(j)}_z}\right)\circ \left(\prod_{i=1}^{m_j}\pi\left(\gamma_{l^{(j)}_i}\right) \right),
$$
the same kind of computations shows that $d_1^*$ is uniformly bounded on $(H_r^2)^p$.
\qed

\smallskip

For $0<r<r_*$, let $\|d_0^*\|_{r}$ and $\|d_1^*\|_{r}$ denote the corresponding operator norms on $(H^2_r)^k$ and $(H^2_r)^p$ respectively. With \rl{actionsurH} and \rl{lemma_d0}, let us define
\begin{equation}\label{Kr_op_norm}
\CK=\CK_{r_*}:= \max\{\|d_0\|_{r_*}, \|d_1\|_{r_*}, \|d_0^*\|_{r_*}, \|d_1^*\|_{r_*}\}.
\end{equation}

\medskip

Recall $S=\{\gamma_1,\cdots,\gamma_k\}$ the set of generators of the group $G$ and their inverses. For $u: S\to H_r^2$, $0<r<r_*$, let us define
$$\CU:=\left(u(\gamma_l)\right)_{1\leq l\leq k}\in ( H_r^2)^k.$$
By the definition of norm (\ref{norm-Gr}), we notice that
$\|\CU\|_r= \|u\|_{S,r}$.

\begin{lemma}\label{lemma_d_1}
If $\|u\|_{S,r}$ is sufficiently small, and $\pi_u:G\to {\rm Diff}^\omega(M_r)$ is a $G-$action by diffeomorphisms with
$\pi_u(\gamma)=\Exp\{u(\gamma)\}\circ\pi(\gamma)$ for $\gamma\in S$, then
$\left\|d_1 {\CU}\right\|_{r} \lesssim  \|u\|^{2}_{S,r}$, where the implicit constant in `` $\lesssim$ " depends only on the manifold $(M_{r_*},\kappa)$ and $S$.
\end{lemma}

\proof Since $u: S\to H_r^2$ with $\|u\|_{S, r}$ sufficiently small, according to Lemma \ref{lem_exp-isometry}, we have that, for every $\gamma\in S$,
\begin{eqnarray}
\pi_u(\gamma)&:=& \Exp \{u(\gamma)\} \circ \pi(\gamma)\label{F_gamma}\\
&=&\pi(\gamma)\circ \Exp \left\{D\pi(\gamma)^{-1}\cdot u(\gamma)\circ \pi(\gamma)\right\}\nonumber\\
&=& \pi(\gamma) + D\pi(\gamma)\cdot\left(D\pi(\gamma)^{-1}\cdot u(\gamma)\circ \pi(\gamma)\right)+h(\gamma)\nonumber\\
&=&\pi(\gamma)+\tilde u(\gamma)+h(\gamma),\nonumber
\end{eqnarray}
where $\tilde u(\gamma):=u(\gamma)\circ\pi(\gamma)$ and $h(\gamma)$ is the sum of all higher-order terms of $u(\gamma)$, with
$$\|h(\gamma)\|_{r}\lesssim \|u(\gamma)\|^{2}_{r}\leq \|u\|^{2}_{S,r}, \quad \forall \ \gamma\in S.$$

Recall $\{\CW_j\}_{1\leq j\leq p}$ the set of all relations of $G$, i.e., each $\CW_j=\gamma_{l_1^{(j)}}\cdots\gamma_{l_{m_j}^{(j)}}$ is a word of length $m_j$ such that
$\gamma_{l_1^{(j)}}\cdots\gamma_{l_{m_j}^{(j)}}=e$.
Since both $\pi$ and $\pi_u$ are $G-$actions by diffeomorphisms, we have
\begin{equation}\label{group_identity}
\pi\left(\gamma_{l_1^{(j)}}\right)\circ \cdots \circ \pi\left(\gamma_{l_{m_j}^{(j)}}\right)=\pi_u\left(\gamma_{l_1^{(j)}}\right)\circ \cdots \circ \pi_u\left(\gamma_{l_{m_j}^{(j)}}\right)={\rm Id}.
\end{equation}
From now on, with $j$ fixed, we omit the superscript ``$(j)$" for convenience. In view of (\ref{F_gamma}), we have that
\begin{eqnarray*}
\pi_u(\gamma_{l_1})\circ \cdots \circ \pi_u(\gamma_{l_{m_j}})&=& \left(\pi(\gamma_{l_1})+\tilde u(\gamma_{l_1})+h(\gamma_{l_1})\right)\circ  \cdots \circ   \left(\pi(\gamma_{l_{m_j}})+\tilde u(\gamma_{l_{m_j}})+h(\gamma_{l_{m_j}})\right)\\
&=& \pi(\gamma_{l_1}\cdots\gamma_{l_{m_j}})+ \tilde u(\gamma_{l_1})\circ \pi\left(\gamma_{l_2}\cdots\gamma_{l_{m_j}}\right)\\
& & + \,  \left(D\pi(\gamma_{l_1})\left[\pi(\gamma_{l_2})\right]\cdot\tilde u(\gamma_{l_2})\right)\circ \pi\left(\gamma_{l_3}\cdots\gamma_{l_{m_j}}\right)\\
& & + \, \left(D\pi(\gamma_{l_1}\gamma_{l_2})\left[\pi(\gamma_{l_3})\right]\cdot \tilde u(\gamma_{l_3})\right)\circ \pi\left(\gamma_{l_4}\cdots\gamma_{l_{m_j}}\right)\\
 & & + \, \cdots +  \left(D\pi(\gamma_{l_1}\gamma_{l_2} \cdots \gamma_{l_{n-1}})\left[\pi(\gamma_{l_n})\right]\cdot  \tilde u(\gamma_{l_n})\right)\circ \pi\left(\gamma_{l_{n+1}}\cdots\gamma_{l_{m_j}}\right)\\
  & & + \, \cdots +  D\pi(\gamma_{l_1}\gamma_{l_2} \cdots \gamma_{l_{m_j-1}})\left[\pi\left(\gamma_{l_{m_j}}\right)\right]\cdot\tilde u\left(\gamma_{l_{m_j}}\right) + \tilde h({\CW}_j),
  \end{eqnarray*}
where $D\pi(\gamma)\left[\pi(\gamma')\right]\cdot\tilde u(\gamma'')$ means the differential of $\pi(\gamma)$ evaluated at $\pi(\gamma')$ and applied on $\tilde u(\gamma'')$, and $\tilde h({\CW}_j)$ is the sum of all higher-order terms of $u$, satisfying
$$\|\tilde h({\CW}_j)\|_{r}\lesssim \|u\|^{2}_{S,r}.$$

Since, for every $\gamma\in S$, $ \tilde u(\gamma)=u(\gamma)\circ \pi(\gamma)$, we see from (\ref{group_identity}) that, for $j=1,\cdots r$,
\begin{eqnarray*}
{\rm Id} &=& {\rm Id}+  u(\gamma_{l_1})\circ \pi(\gamma_{l_1}\cdots\gamma_{l_{m_j}})\\
& & + \, \left(D\pi(\gamma_{l_1})\cdot u(\gamma_{l_2})\right)\circ \pi(\gamma_{l_2}\cdots\gamma_{l_{m_j}})\\
& &+ \, \left(D\pi(\gamma_{l_1}\gamma_{l_2})\cdot u(\gamma_{l_3})\right)\circ \pi(\gamma_{l_3}\cdots\gamma_{l_{m_j}})\\
  & & + \, \cdots +  \left(D\pi(\gamma_{l_1}\gamma_{l_2} \cdots \gamma_{l_{m_j-1}})\cdot u(\gamma_{l_{m_j}})\right)\circ \pi(\gamma_{l_{m_j}}) + \tilde h({\CW}_j)\\
  &=& {\rm Id}+ (d_1{\CU})_j+ \tilde h({\CW}_j),
\end{eqnarray*}
recalling the definition of $d_1$ given in (\ref{defi_d1}).
Then $(d_1{\CU})_j= -\tilde h({\CW}_j)$, and hence $$\|d_1{\CU}\|_{r}\lesssim \|u\|^{2}_{S,r}.\qed$$

\subsection{Kernel of $\square$}

Recall the definition of the box operator $\square$ on $(L^2(M,TM))^k$ in (\ref{def-box}). Then we have the orthogonal decomposition
$$(L^2(M,TM))^k= {\rm Ker} \square \bigoplus {\rm Im} \square.$$

\begin{defi}\label{def-harmonic} The {\bf harmonic operator} $\CH$ on $(L^2(M,TM))^k$ is the projection onto ${\rm Ker} \square$. Let ${\CH}^{\perp}:=\id-\CH$.
\end{defi}

Given $u\in(H_r^2)^k$ with $r\in ]0,r_*[$, according to Remark \ref{rmk_L2}, we have that $$\left. u\right|_{M}\in (\Gamma^{\omega})^k\subset (L^2(M,TM))^k,$$
then ${\CH}u$ and ${\CH}^{\perp}u$ are naturally defined on $M_r$ as the extensions of ${\CH}\left(\left. u\right|_{M}\right)$ and ${\CH}^{\perp}\left(\left. u\right|_{M}\right)$ respectively. 

 Recalling the projection $\PP_i$ defined in (\ref{Proj_i}),
we have
$$\PP_i \circ {\CH} ={\CH} \circ \PP_i ,\qquad \PP_i \circ {\CH}^{\perp}={\CH}^{\perp}\circ \PP_i ,$$
since the operator $\square$ is invariant on $V_i$.


\begin{lemma}\label{lem-normharmonic}
For $r\in ]0,r_*[$ and $u\in(H_r^2)^k$, we have ${\CH}u$,  ${\CH}^{\perp}u\in (H_r^2)^k$, with $\|{\CH}u\|_{r}$, $\|{\CH}^{\perp}u\|_{r}\leq \|u\|_{r}$, and, for every $i\in \N$,
$$\|\PP_i u\|_r\leq e^{r\lambda_i} \|\PP_i  u\|_{L^2},\quad
|\PP_i  u|_{0,r}\leq(1+\lam_i)^{n+1}  e^{r\lambda_i} \|\PP_i  u\|_{L^2}. $$
\end{lemma}


\proof For $1\leq l \leq k$, there exist the coefficients $u_{l,j}$, $j\in I_i$, 
such that
$$\left(\PP_i u \right)_{l} =\sum_{j\in I_i}u_{l,j} {\bf e}_j,$$ and, recalling the definition of $L^2-$norm in \rd{defi_norm_L2pon}, we have
$$\left\|\PP_i u\right\|_{L^2}
=\left(\sum_{1\leq l\leq k}\sum_{j\in I_i}\left|u_{l,j}  \right|^2 \int_{M}  \la {\bf e}_j , {\bf e}_j \ra_g d{\rm vol}  \right)^{\frac12}
=\left(\sum_{1\leq l\leq k}\sum_{j\in I_i}\left|u_{l,j}\right|^2 \right)^{\frac12}. $$
Then, according to the definition of weighted $L^2-$norm in \rd{defi_norm_L2pon}, for any $r\in ]0,r_*[$, we have $\PP_i u\in (H_r^2)^k$ and 
$$\left\|\PP_i u\right\|_{r}= \left(\sum_{1\leq l\leq k}\sum_{j\in I_i}\left|u_{l,j}\right|^2 e^{2r\lambda_i} (1+\lambda_i)^{-\frac{n-1}{2}}\right)^{\frac12}\leq e^{r\lambda_i} \|\PP_i u\|_{L^2}.$$
Furthermore, recalling that $\left|{\bf e}_j\right|_{0,r}\lesssim (1+\lam_i)^{n+1} e^{r\lambda_i} $, we have 
\begin{eqnarray*}
\left|\PP_i u\right|_{0,r}&=&\left(\sum_{1\leq l\leq k}\left|\sum_{j\in I_i}u_{l,j} {\bf e}_j\right|_{0,r}^2 \right)^{\frac12}\\
&\leq&\sup_{j\in I_i}\left|{\bf e}_j\right|_{0,r}\left(\sum_{1\leq l\leq k}\left(\sum_{j\in I_i}\left|u_{l,j}\right|\right)^2 \right)^{\frac12} \ \lesssim \ (1+\lam_i)^{n+1} e^{r\lambda_j} \|\PP_i u\|_{L^2}.
\end{eqnarray*}
On the other hand, for every $i\in \N$, we have
\begin{eqnarray*}
\left\|\PP_i u\right\|_{r}&=&e^{r\lambda_i} (1+\lambda_i)^{-\frac{n-1}{4}} \left(\sum_{1\leq l\leq k} \sum_{j\in I_i}\left|u_{l,j}\right|^2\right)^{\frac12}\\
&=&e^{r\lambda_i} (1+\lambda_i)^{-\frac{n-1}{4}} \|\PP_i u\|_{L^2}\\
&\geq & e^{r\lambda_i} (1+\lambda_i)^{-\frac{n-1}{4}} \|(\CH\circ\PP_i) u\|_{L^2}\\
&\geq&e^{r\lambda_i} (1+\lambda_i)^{-\frac{n-1}{4}} \|(\PP_i \circ {\CH})  u\|_{L^2}
\ = \ \|(\PP_i \circ {\CH})  u\|_{r}.
\end{eqnarray*}
Hence, $\left\|{\CH} u\right\|_{r}\leq  \|u\|_{r}$. Similarly, we have $\left\|{\CH}^{\perp} u\right\|_{r} \leq   \|u\|_{r}$.\qed

\medskip

\begin{lemma}\label{lemmaH_d0} For any $r\in ]0,r_*[$, ${\CH}\circ d_0=0$ on $H^2_r$.
\end{lemma}
\proof
It is sufficient to show that ${\CH}\circ d_0=0$ on $L^2(M,TM)$, and, in view of the Peter-Weyl decomposition (\ref{peter-weyl}), it is sufficient to show that, for $\square_i=\square\circ\PP_i$, $i\in \N$,
$$ d_0v \in{\rm Im}\square_i,\quad \ \forall  \  v \in V_i.$$

For $i\in \N$, since  $V_i^k={\rm Ker}\square_i\bigoplus{\rm Im}\square_i$,
we have the unique decomposition of $d_0v$ for every $v\in V_i$:
$$d_0v=w^{\rm Ker}+w^{\rm Im},\qquad w^{\rm Ker} \in {\rm Ker} \square_i,\quad w^{\rm Im} \in {\rm Im} \square_i,$$
which implies that
$$
\square(d_0v-w^{\rm Im})=\left(\left(d_0\circ d_0^*+d_1^*\circ d_1\right)\circ d_0\right)v-\square w^{\rm Im}
=\left(d_0\circ d_0^*\circ d_0\right)v-\square w^{\rm Im}=0.
$$
Recalling that $v$ is arbitrarily chosen in $V_i$, we have that ${\rm Im} (d_0\circ d_0^*\circ d_0\circ \PP_i)\subset{\rm Im} \square_i.$

Now let us show that $d_0 v\in {\rm Im}\left.(d_0\circ d_0^*\circ d_0)\right|_{V_i}$ for $v\in V_i$, which will complete the proof.
Since $V_i={\rm Ker}(d_0\circ \PP_i)\bigoplus{\rm Im}(d_0^*\circ \PP_i)$, we have the unique decomposition for $v$:
$$v=v^{\rm Ker}+v^{\rm Im},\qquad v^{\rm Ker}\in {\rm Ker}(d_0\circ \PP_i), \quad v^{\rm Im}\in {\rm Im}(d_0^*\circ \PP_i).$$
Then, there exists some $ u\in V_i^k$ satisfying $d_0^*u=v^{\rm Im}$ such that
$$d_0v=d_0v^{\rm Im}=(d_0\circ d_0^*)u.$$
In view of the fact that $V_i^k={\rm Ker}(d_0^*\circ \PP_i)\bigoplus{\rm Im}(d_0\circ \PP_i)$, we have the unique decomposition
$$u=u^{\rm Ker}+u^{\rm Im},\qquad u^{\rm Ker}\in {\rm Ker}(d_0^*\circ \PP_i), \quad u^{\rm Im}\in {\rm Im}(d_0\circ \PP_i).$$
Hence, there exists some $ z\in V_i$ satisfying $d_0 z=u^{\rm Im}$ such that
$$ d_0v=(d_0\circ d_0^*)u=(d_0\circ d_0^*)u^{\rm Im}=(d_0\circ d_0^*\circ d_0)z.$$
This implies that $d_0v\in {\rm Im}(d_0\circ d_0^*\circ d_0\circ \PP_i)\subset{\rm Im} \square_i$.\qed

\subsection{Inverse of $\square$}\label{sec_inverse_box}

\begin{prop}\label{prop_sol_cohomo}
Given $u\in (H_r^2)^k$, $r\in ]0,r_*[$, with ${\CH}u=0$, for any $0<r'<r$, there exists a unique $f\in (H_{r'}^2)^k$ with ${\CH}f=0$ such that $\square f=u$ with
\begin{equation}\label{cohomo_esti}
\|f\|_{r'}\lesssim\frac{\|u\|_{r}}{(r-r')^{\tau}}.
\end{equation}
\end{prop}

\proof For $u\in (H_r^2)^k\subset(L^2(M,TM))^k$, since ${\CH}u=0$,
according to the Peter-Weyl decomposition (\ref{peter-weyl}), it is sufficient to solve the equation
\begin{equation}\label{cohomo_decomp}
\square_i f_i = \left(\sum_{j\in I_i} u_{l,j} {\bf e}_j \right)_{1\leq l\leq k},\qquad i\in \N.
\end{equation}
In view of (\ref{u_EquiForms}), combining with the Diophantine condition (\ref{diophantine}) of $\pi$, Eq. (\ref{cohomo_decomp}) is solvable with the solution
$$\left(\sum_{j\in I_i} f_{l,j} {\bf e}_j\right)_{1\leq l\leq k} = f_i=\sum_{1\leq j\leq K_i}\frac{\tilde u_{i,j}}{\mu_{i,j}} \CE_{i,j}, $$
which implies that, through (\ref{u_EquiForms-norm})
$$\left(\sum_{1\leq l\leq k}\sum_{j\in I_i} \left|f_{l,j} \right|^2 \right)^{\frac12} = \left\|f_i\right\|_{L^2}= \left(\sum_{1\leq j\leq K_i}\left|\frac{\tilde u_{i,j}}{\mu_{i,j}}\right|^2\right)^{\frac12}.$$
Recalling the Diophantine condition (\ref{diophantine}), 
by Definition \ref{defi_norm_L2pon} of the weighted $L^2-$norm, we have that, $f_i\in (V_i\cap H_{r'}^2)^k$ with
\begin{itemize}
\item if $K_0>0$, then
$$\left\|f_0\right\|_{r'}=\left\|f_0\right\|_{L^2}
=\left(\sum_{j=1}^{K_0}\left|\frac{\tilde u_{i,j}}{\mu_{i,j}}\right|^2\right)^{\frac12}
\leq\sigma^{-1}\left(\sum_{j=1}^{K_0}\left|\tilde u_{i,j}\right|^2\right)^{\frac12}
=\sigma^{-1} \left\|u_i\right\|_{r},$$
\item for $i>0$,
\begin{eqnarray*}
\left\|f_i\right\|_{r'}&=&e^{r'\lam_i}(1+\lam_i)^{-\frac{n-1}{4}}\left(\sum_{1\leq l\leq k} \sum_{j\in I_i}|f_{l,j}|^2\right)^{\frac12} \\
&=&e^{r'\lam_i}(1+\lam_i)^{-\frac{n-1}{4}}\left(\sum_{j=1}^{K_i}\left|\frac{\tilde u_{i,j}}{\mu_{i,j}}\right|^2\right)^{\frac12} \\
&\leq&\sigma^{-1}(1+\lam_i)^{\tau}   e^{r'\lam_i}(1+\lam_i)^{-\frac{n-1}{4}} \left(\sum_{j=1}^{K_i}\left|\tilde u_{i,j}\right|^2\right)^{\frac12}\\
&\lesssim&e^{-(r-r')\lam_i}\lam_i^{\tau} \cdot e^{r\lam_i}(1+\lam_i)^{-\frac{n-1}{4}}\left(\sum_{1\leq l\leq k}\sum_{j\in I_i} \left|u_{l,j} \right|^2 \right)^{\frac12}
 \ \lesssim \ \frac{\left\|u_i\right\|_{r}}{(r-r')^{\tau}},
\end{eqnarray*}
since the function $x \mapsto e^{-(r-r') x}x^{\tau} $ attains its maximum on $\R_+$ at $x=\frac{\tau}{r-r'}$, which implies that
$$ e^{-(r-r')\lambda_i}\lambda_i^{\tau} \leq \left(\frac{\tau}{r-r'}\right)^{\tau} e^{-\tau}.$$
\end{itemize}

With the above $f_i$'s, the solution to Eq. (\ref{cohomo_decomp}), being the projection of $f$ onto $(V_i\cap H_{r'}^2)^k$ for every $i\in \N$, and taking $\CH f=0$, the solution $f\in (H_{r'}^2)^k$ to the equation $\square f=u$ is unique. \qed

\medskip

By Proposition \ref{prop_sol_cohomo}, $\square$ is invertible on ${\rm Im}\square$.
Hence, for $0<r'<r<r_*$, the operator
\begin{equation}\label{boxInverse}
\begin{array}{rccl}
\CH^{\perp}\circ  \square^{-1}\circ \CH^{\perp}:&(H_r^2)^k&\to & (H^2_{r'})^k \\
&u&\mapsto& f \  \ {\rm with} \ \CH f=0  \  \  {\rm s. t.} \  \  \square f=\CH^{\perp} u
\end{array}
\end{equation}
is well defined and bounded (depending on $r-r'$ as shown in (\ref{cohomo_esti})). In the sequel, we still use $\square^{-1}$ to denote the above operator for convenience.

\section{KAM scheme -- Proof of Theorem \ref{thm_conj} and \ref{thm-geom}}\label{sec-KAM}

 Let $\pi_0$ be a $G-$action with $\pi_0(\gamma)=\Exp\{P_0(\gamma)\}\circ\pi(\gamma)$ for $\gamma\in S$, where $P_0:S\to H_{r_0}^2$, $r_0\in ]0,r_*[$, with $\|P_0\|_{S,r_0}=\varepsilon_0$ and $\pi$ satisfies the Diophantine condition (\ref{diophantine}). Let ${\CI}:=[\frac{r_0}{2}, r_0]\subset ]0,r_*[$.

In Theorem \ref{thm-geom}, we assume that the first cohomology group of the complex \re{complex} 
is vanishing, which implies, according to \rl{lem-H1}, that $ \Ker\;\square=0$.
 
 In Theorem \ref{thm_conj}, we assume that $\pi_0\in \CF^M_{\pi,\zeta}$ with 
  $0<\zeta<|\ln(\varepsilon_0)|^{-\frac13}$, and $\pi$ satisfies that 
  $\dim\Ker\;\square<\infty$, which means there exists $\CJ\in\N$, such that
 \begin{equation}\label{dimensionKer}
 \Ker\;\square\subset \bigoplus_{i\in \N\atop{\lambda_i\leq \CJ}} V_i^k.
 \end{equation}
By Lemma \ref{lem-normharmonic}, for any $r\in\CI$, $u\in (H_r^2)^k $,
 $$\| {\CH} u\|_r\leq e^{r\CJ} \|{\CH}u\|_{L^2},\quad
|{\CH}u|_{0,r}\leq(1+\CJ)^{n+1}  e^{r\CJ} \|{\CH} u\|_{L^2}. $$
 
We shall prove Theorem \ref{thm_conj} and \ref{thm-geom}  by performing the KAM schemes, which have many common procedures.
 
\subsection{Sequences in KAM scheme}\label{KAM-scheme}

Let $c_{\lesssim}$ be the maximal implicit constant, depending on the manifold $(M_{r_*},\ka)$, the group $G$, the Diophantine constants $\sigma$, $\tau$, and the interval ${\CI}$, appearing in all inequalities involving an inequality of the form ``$\lesssim$".

Let us assume that $\varepsilon_0$ is sufficiently small such that
\begin{equation}\label{vareps0_small}
\left(60+\CJ+e^{r_0\CJ} +\CK +n+ k+ p+ c_{\lesssim} \right)^{12+9n} \left(\frac{24 e^2}{r_0} \right)^{\tau+6n}  <|\ln(\varepsilon_0)|<\varepsilon_0^{-\frac{1}{60(\tau+6n)}},
\end{equation}
with the upper bound $\CK$ of operator norms defined in (\ref{Kr_op_norm}), and, related to the finitely presented group $G$, $k=\# S$, $p=\# R$ given in Section \ref{subsec_GroupAction}.
Note that (\ref{vareps0_small}) implies that 
\begin{equation}\label{use_vareps0}
\left(1+\CJ\right)^{6n+6} e^{6r_0\CJ}<\left(1+\CJ+e^{r_0\CJ} \right)^{12+6n}<|\ln(\varepsilon_0)|,\qquad c_{\lesssim}<|\ln(\varepsilon_0)|^{\frac1{12}}.
\end{equation}


Recalling $\varepsilon_*, \varepsilon_{M}\in ]0,1[$ obtained in Proposition \ref{lemMoser1}, let $\varepsilon_{**}\in ]0,\min\{\varepsilon_*, \varepsilon_{M}\}[$ such that Lemma \ref{lem_exp-isometry} and Lemma \ref{lemma_d_1} are applicable. Assume that $\varepsilon_0$ also satisfies that
\begin{equation}\label{vareps0_small-add}
 \frac{2^{5n+3} |\ln(\varepsilon_0)|^{\frac14}}{r_0^{3n+1}} \cdot \varepsilon_0\leq \varepsilon_{**}.\end{equation}

Let $\zeta_0:=\zeta\in\left]0,|\ln(\varepsilon_0)|^{-\frac13}\right[$. We define the sequences $\{r_m\}$, $\{\tilde r_m\}$, $\{r'_m\}$, $\{\varepsilon_m\}$, $\{N_m\}$, $\{\zeta_m\}$ by
\begin{equation}\label{seqs}
\begin{array}{lll}
\displaystyle r_{m+1}:= r_m- \frac{r_0}{2^{m+2}} , & \displaystyle \tilde r_m:= \frac{r_{m}+r_{m+1}}{2},  & \displaystyle  r_m':=\frac{3r_{m}+r_{m+1}}{4}, \\[3mm]
 \displaystyle \varepsilon_{m+1} := \varepsilon_m^{\frac65} = \varepsilon_0^{(\frac65)^{m+1}}, &\displaystyle  N_m:=\frac{6|\ln(\varepsilon_m)|}{r_m-r_{m+1}},  & \displaystyle \zeta_{m+1}:=\zeta_{m}+2\varepsilon_m^{\frac34}
 \end{array}.
\end{equation}
It is easy to see that,  for every $m\in\N$,
\beq\label{zeta01}
0<  \zeta_m\leq \zeta_0+2\sum_{j\geq 1} \varepsilon_j^{\frac34}<\zeta_0+ \varepsilon_0^{\frac12}<\frac14|\ln(\varepsilon_0)|^{-\frac14},
\end{equation}
$r_m\in \CI=[\frac{r_0}{2}, r_0]$, $r_{m+1}<\tilde r_{m}<r_m'<r_m$ with
\beq\label{dist_rr+}
r_m-r_m'= r_m'- \tilde r_m= \frac{r_{m}-r_{m+1}}{4},\quad \tilde r_m-r_{m+1}=\frac{r_{m}-r_{m+1}}{2},\\
\eeq
and $r_m\to \frac{r_0}{2}$ as $m\to\infty$.



\begin{lemma}\label{prop_seq} For every $m\in\N^*$,
\begin{eqnarray}
\frac{N_m^\tau}{(r_m-r_{m+1})^{6n}} \cdot \varepsilon_m^{\frac{1}{30}}&<&1,\label{vareps-m_small}\\
\frac{2|\ln(\varepsilon_0)|^{\frac14}}{(r_m-r_{m+1})^{3n+1}} \cdot \varepsilon_m &<& \varepsilon_{**}.\label{vareps-m_small-add}
\end{eqnarray}
\end{lemma}
\proof According to the definition of $\{r_m\}$ and $\{N_m\}$ in (\ref{seqs}), we see that
$$\frac{N_m^\tau}{(r_m-r_{m+1})^{6n}}=\frac{6^\tau|\ln(\varepsilon_0)|^\tau}{r_0^{\tau+6n}}
\left(\frac65\right)^{m\tau}2^{(m+2)(\tau+6n)}<\left(\frac{24|\ln(\varepsilon_0)|}{r_0}\right)^{\tau+6n}\left(\frac{12}{5} \right)^{m(\tau+6n)}.$$
Then, under the assumption (\ref{vareps0_small}), we have
$$
\ln\left(\frac{N_m^\tau}{(r_m-r_{m+1})^{6n}} \right)< (\tau+6n)\left(\ln\left(\frac{24}{r_0}\right) + \ln(|\ln(\varepsilon_0)|)  +m \right)< \frac{|\ln(\varepsilon_{0})|}{40}\left(\frac65\right)^m
$$
which implies (\ref{vareps-m_small}).

Since (\ref{vareps-m_small}) implies that
$$\frac{\varepsilon_m}{(r_m-r_{m+1})^{3n+1}}\leq\varepsilon_m^{\frac{29}{30}}\cdot \frac{\varepsilon_m^{\frac{1}{30}}}{(r_m-r_{m+1})^{6n}}\leq  \varepsilon_m^{\frac{29}{30}}\leq   \varepsilon_0^{\frac{29}{25}}\leq  \frac{2^{5n+2} \varepsilon_0}{r_0^{3n+1}},$$
we obtain (\ref{vareps-m_small-add}) through (\ref{vareps0_small-add}).
 \qed

\medskip

We shall present a KAM scheme for the $G-$action $\pi_0$ with 
$$\pi_0(\gamma)=\Exp\{P_0(\gamma)\}\circ \pi(\gamma),\qquad \gamma\in S, $$
under the assumption of Theorem \ref{thm-geom} or \ref{thm_conj}.
According to (\ref{vareps0_small-add}), combining with (\ref{supnormequiv-2}) and Lemma \ref{lem_norm-1}, $\|P_0\|_{S,r_0}<\varepsilon_0$ implies that
\begin{equation}\label{P0-smallnorms}
\|P_0\|_{S,1,r'_0}<\varepsilon_*,\qquad   \|P_0\|_{S,1,M}<\varepsilon_M.\end{equation}
For given $m\in\N$, let us assume that, there is $P_m:S\to H^2_{r_m}$ satisfying
\begin{equation}\label{hypo_Pm}
\|P_m\|_{S,r_m}=\left(\sum_{1\leq l\leq k}\|P_m(\gamma_l)\|_{r_m}^2\right)^{\frac12}<\varepsilon_m,
\end{equation}
such that, $\pi_0$ is conjugate to some $G-$action $\pi_m$, which can be written as 
$$\pi_m(\gamma)=\Exp\{P_m(\gamma)\}\circ \pi(\gamma),\qquad \gamma\in S. $$ In view of (\ref{vareps-m_small-add}), combining with (\ref{supnormequiv-2}) and Lemma \ref{lem_norm-1}, we have that
\begin{equation}\label{Pm-smallnorms}
\|P_m\|_{S,1,r'_m}<\varepsilon_*,\qquad  \|P_m\|_{S,1,M}<\varepsilon_M.
\end{equation}
We shall show that, by a biholomorphism $\Exp\{w_m\}$ on $M_{r_{m+1}}$,
for every $\gamma\in S$, $\pi_m$ can be conjugated to some $G-$action $\pi_{m+1}$, written as 
$$\pi_{m+1}(\gamma)=\Exp\{P_{m+1}(\gamma)\}\circ \pi(\gamma), \qquad \gamma\in S,$$ with $P_{m+1}:S\to H^2_{r_{m+1}}$ satisfying that
$\|P_{m+1}\|_{S,r_{m+1}}<\varepsilon_{m+1}$.
Then Theorem \ref{thm_conj} and \ref{thm-geom}  will be shown by the convergence of the above iterative procedure.

\subsection{Truncation operator} As a first procedure of one KAM step, we define the truncation operator on the Hardy space and introduce several properties.
\begin{defi}\label{def_truncation}
Given $\nu , N\in\N^*$, let us define the {\bf truncation operator} of degree N, ${\Pi}^{(\nu)}_{N}: \left(L^2(M,TM)\right)^\nu\to \left(L^2(M,TM)\right)^\nu$ as
$$
{\Pi}^{(\nu)}_{N}u:=\bigoplus_{i\in \N\atop{\lambda_i\leq N}}\PP_i u=\left(\sum_{i\in \N\atop{\lambda_i\leq N}}\sum_{j\in I_i} u_{l,j} {\bf e}_j\right)_{1\leq l\leq \nu}, \quad u=\left(\sum_{j\in\N} u_{l,j}{\bf e}_j\right)_{1\leq l\leq \nu} \in \left(L^2(M,TM)\right)^\nu,
$$
and let ${\Pi}^{(\nu)\perp}_{N}:=\id-{\Pi}^{(\nu)}_{N}=\bigoplus_{i\in \N\atop{\lambda_i> N}}\PP_i$.
\end{defi}
\noindent It is easy to see that
\begin{equation}\label{ImTr-V_N}
{\rm Im}\Pi_N^{(\nu)}= \bigoplus_{i\in \N\atop{\lam_i\leq N}} V^\nu_i =:\V^\nu_{N}.
\end{equation}
As the harmonic operator $\CH$ (recall Definition \ref{lem-normharmonic}), for given $r\in ]0,r_0[$ and $u\in(H_r^2)^k$, ${\Pi}^{(\nu)}_N u$ and ${\Pi}^{(\nu)\perp}_N u$ are naturally defined by identifying the coefficients with those of the restrictions on $M$.

\smallskip

With $N_m$ defined in (\ref{seqs}), let us now focus on the truncation operator ${\Pi}^{(k)}_{N_m}$, as well as ${\Pi}^{(k)\perp}_{N_m}$, on $(H^2_{r_m})^{k}$.
For convenience, the superscript ``$(k)$" in the notation of the truncation operator will be omitted in the sequel since there is no ambiguity.

Let
${\CP}_m:=\left(P_m(\gamma_l)\right)_{1\leq l\leq k}\in (H_{r_m}^2)^k$.
\begin{lemma}\label{lemma_error_trunc}
${\Pi}^{\perp}_{N_m}\circ\CH=0$  and  $\left\|{\Pi}^{\perp}_{N_m}\CP_m\right\|_{r'_{m}}< \varepsilon_{m}^{2}$.
\end{lemma}
\proof In view of (\ref{dimensionKer}), we have that $\PP_i\circ\CH=0$ whenever $\lambda_i\geq \CJ+1$. Since the inequality (\ref{vareps0_small}) implies that the sequence $\{N_m\}$ defined in (\ref{seqs}) satisfies 
$$N_m\geq N_0 =\frac{6|\ln(\varepsilon_0)|}{r_0-r_{1}}= \frac{24|\ln(\varepsilon_0)|}{r_0}\geq \CJ+1, \qquad \forall \ m\in\N,$$
we have ${\Pi}^{\perp}_{N_m}\circ\CH=\bigoplus_{i\in \N\atop{\lambda_i> N_m}}(\PP_i\circ\CH)=0$.

According to the assumption \re{hypo_Pm}, for $1\leq l\leq k$, $P_m(\gamma_l)=\sum_{j\in\N} P_{m,j}(\gamma_l) {\bf e}_j$ satisfies that $\|P_m(\gamma_l)\|_{r_{m}}<\varepsilon_m$. Then, by (\ref{decay_coeff}) in Proposition \ref{prop_estim-norm}, for every $i\in \N$,
$$|P_{m,j}(\gamma_l)|\leq \varepsilon_m e^{-r_{m}\lambda_i} (1+\lambda_i)^{\frac{n-1}{4}},\quad j\in I_i.$$
In view of Definition \ref{def_truncation}, we see that
$${\Pi}^{\perp}_{N_m}\CP_m=\left(\sum_{i\in \N\atop{\lambda_i> N_m}}\sum_{j\in I_i}  P_{m,j}(\gamma_l) {\bf e}_j\right)_{1\leq l\leq k}.$$
By (\ref{asymp-lamdak}), we have that, for every $K\in\N^*$,
$$\# \bigcup_{i\in \N\atop{\lambda_i\leq K}}I_i\lesssim K^n .$$
Recalling \re{dist_rr+}, we obtain that
\begin{eqnarray*}
\left\|{\Pi}^{\perp}_{N_m}\CP_m\right\|^2_{r'_{m}} &=&\sum_{1\leq l\leq k} \sum_{i\in \N\atop{\lambda_i> N_m}}\sum_{j\in I_i} |P_{m,j}(\gamma_l)|^2 e^{2r'_{m}\lambda_i} (1+\lambda_i)^{-\frac{n-1}{2}}\\
&\lesssim&  \varepsilon^2_m \sum_{K= N_m}^\infty (K+1)^n \sum_{K< \lambda_i\leq K+1}e^{-\frac12(r_{m}-r_{m+1})\lambda_i}  \\
&\lesssim& \varepsilon^2_m \int_{N_m}^{+\infty} (t+1)^n e^{-\frac12(r_{m}-r_{m+1})t} \, dt.
\end{eqnarray*}
By multiple integration by part, we have
\begin{eqnarray*}
 \int_{N_m}^{+\infty} (t+1)^n e^{-\frac12(r_{m}-r_{m+1})t} \, dt&\lesssim& e^{-\frac{N_{m}}2(r_{m}-r_{m+1})}\sum_{j=0}^{n}\frac{n! }{j!}\frac{(N_m+1)^j}{(r_{m}-r_{m+1})^{n+1-j}}  \\
&\lesssim&e^{-\frac{N_{m}}2(r_{m}-r_{m+1})}\left(N_m+1+\frac{1}{r_{m}-r_{m+1}}\right)^{n+1} \ \leq \ \varepsilon_m^{\frac52}.
\end{eqnarray*}
The last inequality follows from the fact that
$$e^{-\frac{N_{m}}2(r_{m}-r_{m+1})}=\exp\{-3|\ln(\varepsilon_m)|\}=\varepsilon_m^3$$
and that Lemma \ref{prop_seq} implies
$$\left(N_m+1+\frac{1}{r_{m}-r_{m+1}}\right)^{n+1}\leq  \varepsilon_m^{-\frac12}.$$
This completes the proof through Lemma \ref{prop_seq}.\qed

\medskip

\begin{lemma}\label{lemma_box_tr} For any $r\in ]0,r_*[$, $\square^{-1}\circ{\Pi}_{N_m}:(H_{r}^2)^k\to (H_{r}^2)^k$ is bounded with
\begin{equation}\label{boxT}
\left\|\square^{-1}\circ{\Pi}_{N_m}\right\|_{r}\lesssim  N_m^{\tau}.
\end{equation}
\end{lemma}

\proof Recall $\square^{-1}$, defined in (\ref{boxInverse}), is extended by $0$ on its kernel. 
 According to (\ref{u_EquiForms-norm}), the Diophantine condition (\ref{diophantine}) implies that, for any $u_i=\left(\sum_{j\in I_i} u_{l,j}{\bf e}_j\right)_{1\leq l\leq k}\in (V_i\cap H_{r}^2)^k$,
\begin{eqnarray*}
\left\|\square^{-1} u_i\right\|_{r}&\lesssim&e^{r\lambda_i}(1+\lambda_i)^{-\frac{n-1}{4}} \left( \sum_{1\leq l\leq k}\sum_{1\leq j\leq K_i} \left|\frac{\tilde u_{i,j}}{\mu_{i,j}}\right|^2 \right)^{\frac12}\\
&\lesssim& \frac{(1+\lambda_i)^{\tau}}{\sigma}\left( \sum_{1\leq l\leq k}\sum_{j\in I_i} \left|u_{i,j}\right|^2 \right)^{\frac12}e^{r\lambda_i}(1+\lambda_i)^{-\frac{n-1}{4}}.
\end{eqnarray*}
Summing over $i\in \N$ with $0<\lambda_i\leq N_m$, we obtain  (\ref{boxT}).\qed

\subsection{Cohomological equation -- one KAM step}\label{sect-cohom}

Let us write the equation of the $(m+1)-{\rm th}$ conjugacy for every $\gamma\in S$, i.e., the conjugacy from $\pi_m(\gamma)=\Exp\{P_m(\gamma)\}\circ \pi(\gamma)$ to $\pi_{m+1}(\gamma)=\Exp\{P_{m+1}(\gamma)\}\circ \pi(\gamma)$, by a diffeomorphism $\Exp\{w_m\}$ of $M_{\tilde r_{m}}$, as the following:
\begin{equation}\label{applem-3}
\Exp\{w_m\}^{-1}\circ\Exp\{P_m(\gamma)\}\circ \pi(\gamma) \circ \Exp\{w_m\} = \Exp\{P_{m+1}(\gamma)\}\circ \pi(\gamma),\quad \forall \ \gamma\in S,
\end{equation}
with $w_m\in H^2_{r_m}$, independent of $\gamma$, and $P_{m+1}:S\to H^2_{r'_{m}}$,
 satisfying
\beq\label{estim2reach}
\|w_m\|_{r_m}\leq N_m^{\tau}\varepsilon_{m},\qquad  \|P_{m+1}\|_{S,r'_{m}}\leq \varepsilon_{m}^{\frac54},
\eeq
to be determined.
Let us rewrite this conjugacy equation in a more tractable way. To do so, let us assume that $\|w_m\|_{r_m}$ is small enough 
so that, according to Lemma \ref{lem_exp-isometry}, Eq. \re{applem-3} reads, for every $\gamma\in S$,
\begin{eqnarray*}
\Exp\{P_{m+1}(\gamma)\}&=& \Exp\{w_m\}^{-1} \circ \Exp\{P_m(\gamma)\}\circ \pi(\gamma) \circ \Exp\{w_m\}\circ \pi^{-1}(\gamma)\\
&=& \Exp\{w_m\}^{-1} \circ \Exp\{P_m(\gamma)\} \circ  \Exp\left\{D\pi(\gamma)\cdot w_m\circ\pi^{-1}(\gamma)\right\},
\end{eqnarray*}
and hence
\begin{equation}\label{applem-1}
\Exp\{P_{m}(\gamma)\}\circ  \Exp\left\{D\pi(\gamma)\cdot w_m\circ\pi^{-1}(\gamma)\right\}= \Exp\{w_m\}\circ \Exp\{P_{m+1}(\gamma)\}.
\end{equation}
In view of (\ref{supnormequiv-2}) and Lemma \ref{lem_norm-1}, we see that $\|w_m\|_{1,r_m'}$, $\|P_m\|_{S,1,r_m'}$ and $|P_{m+1}|_{S,0,r_m'}$ would be smaller than $\varepsilon_*$, so that we could apply Proposition \ref{lemMoser1} to both sides of (\ref{applem-1}),
we have, on $M_{\tilde r_{m}}$, for every $\gamma\in S$,
\begin{eqnarray}
   & & P_m(\gamma)+D\pi(\gamma)\cdot w_m\circ\pi^{-1}(\gamma)+s_1\left(P_{m}(\gamma), D\pi(\gamma)\cdot w_m\circ\pi^{-1}(\gamma)\right) \label{eq_cohomo-im}\\
   &= & w_m+P_{m+1}(\gamma)+ s_1\left(w_m, P_{m+1}(\gamma)\right), \nonumber
\end{eqnarray}
where, for every $\gamma\in S$,
$$s_1\left(P_{m}(\gamma), D\pi(\gamma)\cdot w_m\circ\pi^{-1}(\gamma)\right), \ s_1\left(w_m, P_{m+1}(\gamma)\right)\in \Gamma_{\tilde r_m}. $$
and, in view of (\ref{vareps0_small}), for $1\leq l\leq k$, we would have
\begin{eqnarray}
|s_1(w_m, P_{m+1}(\gamma_l))|_{0,\tilde r_m}&\lesssim& |w_m|_{1,r'_{m}}|P_{m+1}(\gamma_l)|_{0,\tilde r_m}, \label{error_s1-1}\\
|s_1(P_{m}(\gamma_l), D\pi(\gamma_l)\cdot  w_m\circ\pi^{-1}(\gamma_l))|_{0,\tilde r_m} &\lesssim& |P_{m}(\gamma_l)|_{1,r'_{m}}|w_m|_{0,\tilde r_m}. \label{error_s1-2}
\end{eqnarray}


Let us set\begin{eqnarray*}
{\CP}_{m+1} &:=&  \left(P_{m+1}(\gamma_l)\right)_{1\leq l\leq k},\\
\CR_{m,1}&:=& \left(s_1(w_m, P_{m+1}(\gamma_l))\right)_{1\leq l\leq k},\\
\CR_{m,2}&:=& \left(s_1(P_{m}(\gamma_l), \pi(\gamma_l)_*w_m )\right)_{1\leq l\leq k}.
\end{eqnarray*}
Recalling that
$d_0 w_m=\left(w_m-\pi(\gamma_l)_*w_m \right)_{1\leq l\leq k},$
Eq. \re{eq_cohomo-im} with $\gamma=\gamma_1,\cdots,\gamma_k\in S$ can be written as the {\bf cohomological equation} on $M_{\tilde r_m}$:
\begin{equation}\label{cohomo}
\CP_m-d_0 w_m=\CP_{m+1}+\CR_{m,1} - \CR_{m,2}.
\end{equation}
Our goal is to find $w_m\in H_{r_{m}}^2$ and $\CP_{m+1}\in (H_{\tilde r_{m}}^2)^k$ satisfying this (implicit) equation in $(H_{\tilde r_{m}}^2)^k$ as well as the aforementioned estimates, then we also solve Eq. \re{applem-3}. To do so, we proceed as follow.

In view of Lemma \ref{lemma_box_tr}, let us set
\begin{equation}\label{f_w}
\CF_m:=\left(\square^{-1}\circ {\Pi}_{N_m}\circ\CH^{\perp}\right)\CP_m \in \left(H_{r_m}^2\right)^k.
\end{equation}
In the sequel, when there is no possible confusion, we shall often write $\square^{-1}$ instead of $\CH^{\perp} \circ \square^{-1}\circ\CH^{\perp}$, following the notation \re{boxInverse}.
By definition, we have that $\CH \CF_m=0$ and $\square \CF_m= \left({\Pi}_{N_m}\circ {\CH}^{\perp}\right)\CP_m$. Let us define
\beq\label{changem} w_m:=d_0^* \CF_m \in H^2_{r_m}.
\eeq
Since the restriction to $M$ of $\CF_m$ is an analytic vector field on $M$, according to formula \re{d_0*v}, so is $w_m$.
In view of (\ref{vareps0_small}), (\ref{boxT}), Lemma \ref{prop_seq}, and by the uniform boundness \re{boundd0*} of $d_0^*$, we obtain
\begin{equation}\label{esti_wm}
\|w_m\|_{r_m} =\|d^*_0\CF_m\|_{r_m}\lesssim N_m^{\tau} \varepsilon_m.
\end{equation}
This implies, through (\ref{supnormequiv-2}) and Lemma \ref{lem_norm-1}, that
\begin{equation}\label{esti-wm-01}
|w_m|_{0,r'_m}\lesssim \frac{\|w_m\|_{r_m}}{(r_m-r_{m+1})^{3n}}\leq \varepsilon_m^{\frac{19}{20}},\quad  |w_m|_{1, r'_m}\lesssim \frac{\|w_m\|_{r_m}}{(r_m-r_{m+1})^{3n+1}}\leq \varepsilon_m^{\frac{19}{20}}.
\end{equation}

Recalling that $\square=d_0\circ d_0^*+d_1^*\circ d_1 $ and \re{f_w}, we have
\begin{eqnarray}
\CP_m-d_0w_m&=&\CP_m- (d_0\circ d_0^*)\CF_m\label{P-d0}\\
   &=&{\CH}\CP_m+{\CH}^{\perp}\CP_m - \square \CF_m+(d_1^*\circ d_1) \CF_m \nonumber\\
   &=&{\CH}\CP_m+{\CH}^{\perp}\CP_m-\left({\Pi}_{N_m}\circ{\CH}^{\perp}\right)\CP_m+(d_1^*\circ d_1) \CF_m \nonumber\\
   &=&{\CH}\CP_m+\left({\Pi}^{\perp}_{N_m}\circ\CH^{\perp}\right)\CP_m +(d_1^*\circ d_1) \CF_m. \nonumber
\end{eqnarray}

\begin{lemma}\label{lemma_d_1d_1} $(d_1^*\circ d_1) \CF_m\in \left( {\Bbb V}_{N_m}\cap H_{r_m}^2\right)^k$ with
 $\|(d_1^*\circ d_1) \CF_m\|_{r_m}\lesssim N_m^{\tau}\varepsilon_{m}^{\frac32}$.
\end{lemma}
\proof Recalling (\ref{f_w}), we see that
$$\CF_m=(\square^{-1}\circ {\Pi}_{N_m})\CP_m\in \left( {\Bbb V}_{N_m}\cap H_{r_m}^2\right)^k.$$
In view of the assertion (i) in Proposition \ref{prop_estim-norm}, we have that
\begin{eqnarray*}
& &\left|\left\la\left(d_1^*\circ d_1 \circ \square^{-1}  \circ {\Pi}_{N_m}\right) \CP_m,   {\Pi}_{N_m} \CP_m \right\ra_{H^2_{r_m}}\right|\\
&=& \left| \left\la \left(d_1 \circ \square^{-1}  \circ {\Pi}_{N_m}\right) \CP_m,   \left(d_1\circ {\Pi}_{N_m}\right)\CP_m   \right\ra_{H^2_{r_m}}\right|\\
&=& \left| \left\la  \left(d_1 \circ \square^{-1}  \circ {\Pi}_{N_m}\right) \CP_m,   \left({\Pi}_{N_m} \circ d_1 \right)\CP_m \right\ra_{H^2_{r_m}}\right|\\
&\lesssim& \left\|\left(d_1 \circ \square^{-1}  \circ {\Pi}_{N_m}\right) \CP_m\right\|_{r_m} \left\| \left({\Pi}_{N_m} \circ d_1 \right)\CP_m \right\|_{r_m}.
\end{eqnarray*}
Through Lemma \ref{lemma_d_1} and \ref{lemma_box_tr}, we obtain that
$$ \left\|\left(d_1 \circ \square^{-1}  \circ {\Pi}_{N_m}\right) \CP_m\right\|_{r_m}\lesssim N_{m}^\tau \|\CP_m\|_{r_m} , \qquad  \|d_1 \CP_m \|_{r_m}\lesssim \|\CP_m\|_{r_m}^{2}$$
Then, we obtain that
\begin{eqnarray*}
& &\left\| \left((d_1^*\circ d_1\circ \square^{-1} )^\frac12\circ {\Pi}_{N_m}\right) \CP_m\right\|_{r_m}\\
&\lesssim& \left\| \left((d_1^*\circ d_1\circ \square^{-1} )^\frac12\circ {\Pi}_{N_m}\right) \CP_m\right\|_{(H^2_{r_m})^k} \\
&=& \left|\left\la \left((d_1^*\circ d_1\circ \square^{-1} )^\frac12\circ {\Pi}_{N_m}\right) \CP_m, \left((d_1^*\circ d_1\circ \square^{-1} )^\frac12\circ {\Pi}_{N_m}\right) \CP_m \right\ra_{(H^2_{r_m})^k}\right|^{\frac12} \\
  &=& \left|\left\la \left(d_1^*\circ d_1\circ \square^{-1} \circ  {\Pi}_{N_m}\right) \CP_m, {\Pi}_{N_m} \CP_m  \right\ra_{(H^2_{r_m})^k}\right|^{\frac12}\\
  &\lesssim& N_{m}^{\frac{\tau}2}\|\CP_m \|^{\frac32}_{r_m}.
\end{eqnarray*}
Hence, we have
\begin{eqnarray*}
\left\| \left(d_1^*\circ d_1\right) \CF_m\right\|_{r_m}&=& \left\| \left(d_1^*\circ d_1\circ \square^{-1}\circ  {\Pi}_{N_m}\right) \CP_m\right\|_{r_m}\\
&\leq& \left\| \left.\left(d_1^*\circ d_1\circ  \square^{-1}\right)^\frac12\right|_{(\V_{N_m}\cap H_{r_m}^2)^k}\right\|_{r_m} \left\| \left((d_1^*\circ d_1\circ \square^{-1} )^\frac12\circ {\Pi}_{N_m}\right) \CP_m\right\|_{r_m}\\
&\lesssim& N_{m}^{\frac{\tau}2}\cdot   N_{m}^{\frac{\tau}2} \|\CP_m \|^{\frac32}_{r_m},
\end{eqnarray*}
which shows the lemma by the uniform boundedness of $d_1^*$.
\qed

\medskip

According to (\ref{supnormequiv-1}) and (\ref{dist_rr+}), Lemma \ref{lemma_error_trunc} and Lemma \ref{lemma_d_1d_1} imply that
\begin{equation}\label{error_TddF}
|{\Pi}^{\perp}_{N_m} \CP_m|_{0,\tilde r_m} \lesssim \frac{\varepsilon_m^2}{(r_m-r_{m+1})^{3n}},\qquad |(d_1^*\circ d_1) \CF_m|_{0,\tilde r_m}\lesssim \frac{N_m^{\tau}\varepsilon_m^{\frac32}}{(r_m-r_{m+1})^{3n}}.
\end{equation}
both of which are bounded from above by $\varepsilon_{m}^{\frac43}$ through Lemma \ref{prop_seq}.
Hence, in view of (\ref{vareps0_small}), (\ref{P-d0}) and Lemma \ref{lem-normharmonic}, we have
\begin{equation}\label{Pd0w}
|\CP_m-d_0w_m|_{0,\tilde r_m}\leq |{\CH}\CP_m|_{0,\tilde r_m} + |{\Pi}^{\perp}_{N_m} \CP_m|_{0,\tilde r_m}+|(d_1^*\circ d_1) \CF_m|_{0,\tilde r_m}
\end{equation}
Recalling (\ref{error_s1-1}), (\ref{error_s1-2}) and (\ref{esti-wm-01}), we have
\begin{equation}\label{esti_s1-ab}
|{\CR}_{m,1}|_{0,\tilde r_{m}}\lesssim \varepsilon_m^{\frac{19}{20}}\cdot |\CP_{m+1}|_{0,\tilde r_m},\qquad
 |{\CR}_{m,2}|_{0,\tilde r_{m}}\lesssim \frac{\varepsilon_m}{(r_m-r_{m+1})^{3n+1}} \cdot  \varepsilon_m^{\frac{19}{20}} < \varepsilon_m^{\frac32}.
\end{equation}
By the fixed point theorem (see e.g., \cite{dieudonne1}[10.1.1]), together with \rp{lemMoser1}-(i), we obtain the existence of solution $\CP_{m+1}\in  (H_{\tilde r_{m}}^2)^k$ of Eq. (\ref{cohomo}). Its restriction to $M$ is also tangent to $M$. Together with (\ref{Pd0w}) and (\ref{esti_s1-ab}), we obtain its estimate~:
\begin{eqnarray}
\left(1-\varepsilon_m^{\frac{19}{20}}\right) |\CP_{m+1}|_{0,\tilde r_{m}}&\lesssim &|\CP_{m+1}+{\CR}_{m,1}|_{0,\tilde r_m}  \nonumber\\
  &=& |\CP_m-d_0 w_m+{\CR}_{m,2}|_{0,\tilde r_m}\nonumber\\
  &\lesssim&   |{\CH}\CP_m|_{0,\tilde r_m} +\varepsilon_{m}^{\frac43}.\label{rough-harmonique}
\end{eqnarray}

 Recalling that the aim of this $m^{\rm th}$ KAM step is to conjugate the $G-$action $\pi_m(\gamma)$ to $\pi_{m+1}(\gamma)=\Exp{P_{m+1}(\gamma)}\circ\pi(\gamma)$ with
$\|P_{m+1}\|_{S,r_{m+1}}<\varepsilon_{m+1}$. It is sufficient to show that
\begin{equation}\label{Pm+10r}
|\CP_{m+1}|_{0,\tilde r_{m}}=\left(\sum_{1\leq l \leq k}|P_{m+1}(\gamma_l)|^2_{0,\tilde r_{m}}\right)^{\frac12}<\varepsilon^{\frac54}_{m}.
\end{equation}
Indeed, (\ref{Pm+10r}) implies, through (\ref{supnormequiv-2}), that
$$\|P_{m+1}(\gamma_l)\|_{r_{m+1}}\lesssim |P_{m+1}(\gamma_l)|_{0,\tilde r_{m}}<\varepsilon^{\frac54}_{m},\quad  1\leq l \leq k,$$
then, according to the definition of $\|\cdot\|_{S,r_{m+1}}-$norm,
$$\|P_{m+1}\|_{S,r_{m+1}}=\left(\sum_{1\leq l \leq k}\|P_{m+1}(\gamma_l)\|^2_{r_{m+1}}\right)^{\frac12}  \leq k |\ln(\varepsilon_0)|^{\frac14} \varepsilon^{\frac54}_{m} <\varepsilon^{\frac65}_{m}=\varepsilon_{m+1},$$
where, in view of \re{use_vareps0}, the implicit constants in the above inequalities with ``$\lesssim$" are all bounded by $|\ln(\varepsilon_0)|^{\frac1{12}}$.

Under the assumption that 
$H^1(G,V_i)=0$, $\forall \  i\in\N$,
we obtain, through Lemma \ref{lem-H1} that, ${\CH}\CP_m=0$. According to \re{rough-harmonique} and since $0< \varepsilon_m^{\frac{19}{20}}<\frac12$, we obtain $|\CP_{m+1}|_{0,\tilde r_{m}}\lesssim \varepsilon_{m}^{\frac43}$.
Therefore, one step of the KAM scheme for Theorem \ref{thm-geom} is done. 

If the first cohomology group does not vanish and since $0< \varepsilon_m^{\frac{19}{20}}<\frac12$, we barely obtain 
 through Lemma \ref{lem-normharmonic}, the rough estimate of $|\CP_{m+1}|_{0,\tilde r_{m}}:$
\begin{equation}\label{esti_rough_Pm+1}
	|\CP_{m+1}|_{0,\tilde r_{m}}\lesssim  |{\CH}\CP_m|_{0,\tilde r_m} +\varepsilon_{m}^{\frac43}\lesssim \|\CP_m\|_{r_{m}}+\varepsilon_{m}^{\frac43} \lesssim \varepsilon_{m}.
\end{equation}
In what follows, we prove that our ``formal rigidity" assumption allows the refine estimate of $\CP_{m+1}$.

\subsubsection{Refined estimate of $\CP_{m+1}$}\label{refined}

Let us turn to the KAM scheme for Theorem \ref{thm_conj}, for which the first cohomology is not vanishing. 
With $(P_{m+1}(\gamma_l))_{1\leq l\leq k}=\CP_{m+1}$, we obtain the $G-$action $\pi_{m+1}$ generated by $\Exp\{P_{m+1}(\gamma)\}\circ\pi(\gamma)$, $\gamma\in S$. We also found there that $w_m$ satisfied the first estimate of \re{estim2reach} while $\CP_{m+1}$ satisfied the coarser estimate \re{esti_rough_Pm+1}. 
We shall refine the ``rough" estimate (\ref{esti_rough_Pm+1}) to (\ref{Pm+10r}) by taking advantage of the formal conjugacy of $\pi_0$ and hence $\pi_m$.



\begin{lemma}\label{lemma_for_conj} If $\pi_m \in \CF_{\pi,\zeta_m}^M$, then $\pi_{m+1}\in \CF^M_{\pi,\zeta_{m+1}}$. 
\end{lemma}
\proof Since $\pi_m \in \CF^M_{\pi,\zeta_{m}}$ and $\pi_m(\gamma) =  \Exp\{P_{m}(\gamma)\}\circ \pi(\gamma)$ for $\gamma\in S$, by Definition \ref{defi_formal_conj}, for any $j\in\N^*$, there exists $y_{m,j}\in \Ga^\om$ such that
\begin{equation}\label{y_form_conj_m}
\|y_{m,j}\|_{1,M}=|y_{m,j}|_{0,M}+|y_{m,j}|_{1,M}<\zeta_{m},
\end{equation}
and, for every $\gamma\in S$,
\begin{equation}\label{eq_form_conj_m}
\Exp\{y_{m,j}\}^{-1}\circ \Exp\{P_{m}(\gamma)\}\circ \pi(\gamma)\circ \Exp\{y_{m,j}\}=\Exp\{z_{m,j}(\gamma)\}\circ \pi(\gamma),
\end{equation}
where $z_{m,j}:S\to  \Ga^\om$ satisfies that
$$\|z_{m,j}\|_{S,1,M}=\left(\sum_{1\leq l\leq k}\|z_{m,j}(\gamma_l)\|^2_{1,M}\right)^{\frac12}<\varepsilon_{m}^j.$$
Combining with (\ref{applem-3}), we have, for every $\gamma\in S$,
\begin{eqnarray*}
& &\Exp\{y_{m,j}\}^{-1} \circ \Exp\{w_m\}^{-1} \circ \Exp\{P_{m+1}(\gamma)\}\circ \pi(\gamma) \circ \Exp\{w_m\}\circ \Exp\{y_{m,j}\}\\
&=& \Exp\{z_{m,j}(\gamma)\}\circ \pi(\gamma).
\end{eqnarray*}
Then, according to Proposition \ref{lemMoser1}, for
$\tilde y_{m+1,j}:=w_m+ y_{m,j} + s_1(w_m, y_{m,j}) \in \Ga^\om$,
we have that
$
 \Exp\{w_m\}\circ \Exp\{y_{m,j}\} =\Exp\{\tilde y_{m+1,j}\}
$,
and hence, for every $\gamma\in S$,
$$\Exp\{\tilde y_{m+1,j}\}^{-1} \circ \Exp\{P_{m+1}(\gamma)\}\circ \pi(\gamma) \circ \Exp\{\tilde y_{m+1,j}\}=\Exp\{z_{m,j}(\gamma)\}\circ \pi(\gamma). $$
Recalling (\ref{esti-wm-01}), we have, through Proposition \ref{lemMoser1}-(ii), Lemma \ref{lem_norm-1} and \re{esti_wm},
$$\|s_1(w_m, y_{m,j})\|_{1,M}\lesssim \|w_m\|_{1,M}(1+\|y_{m,j}\|_{1,M})\lesssim r_m^{-\left(3n+1\right)}N_m^\tau \varepsilon_m \lesssim N_m^\tau \varepsilon_m,$$
since the sequence $\{r_m\}$ is uniformly bounded from below and $M\subset M_{r_m}$.
In view of \re{use_vareps0}, the implicit constants in the above inequalities with ``$\lesssim$" are all bounded by $|\ln(\varepsilon_0)|^{\frac1{12}}$, we have that
\begin{equation}\label{ineq-lesssim}
\|s_1(w_m, y_{m,j})\|_{1,M}\leq |\ln(\varepsilon_0)| N_m^\tau \varepsilon_m \leq \varepsilon^{\frac34}_m,
\end{equation}
Collecting the two last inequalities, we obtain, for every $j\in\N^*$,
$$ \|\tilde y_{m+1,j}\|_{1,M}\leq \|y_{m,j}\|_{1,M} + \|w_m\|_{1,M} +\|s_1(w_m, y_{m,j})\|_{1,M}
\leq \zeta_{m}+\varepsilon^{\frac34}_m\leq \zeta_{m+1}.$$
Hence, $\pi_{m+1} \in\CF_{\pi,\zeta_{m+1}}$.
\qed

\medskip

\begin{prop} If $\pi_m \in \CF_{\pi,\zeta_m}^M$ and $\dim\Ker\square<\infty$, then
$|\CP_{m+1}|_{0,\tilde r_{m}}<\varepsilon_{m}^{\frac54}$. 
\end{prop}
\proof 

In view of (\ref{cohomo}) and (\ref{P-d0}), we see that $\CP_{m+1}={\CH}\CP_m+\tilde \CP_{m+1}$ with
$$\tilde \CP_{m+1}:=\left({\Pi}^{\perp}_{N_m}\circ \CH^{\perp}\right) \CP_m +(d_1^*\circ d_1) \CF_m-\CR_{m,1}+\CR_{m,2}.$$
By the rough estimate $|\CP_{m+1}|_{0,\tilde r_{m}}\lesssim \varepsilon_m$, we have, through (\ref{esti_s1-ab}), that
$$|\CR_{m,1}|_{0,\tilde r_{m}}\lesssim \varepsilon_m^{\frac{19}{20}}|\CP_{m+1}|_{0,\tilde r_{m}}\lesssim \varepsilon_m^{\frac{39}{20}}.$$
Then, combining with \re{error_TddF}, (\ref{esti_s1-ab}) and Remark \ref{rmk_L2}, we obtain that
\begin{equation}\label{Pm+1L2}
\|\tilde \CP_{m+1}\|_{L^2}\lesssim |\tilde \CP_{m+1}|_{0,M}\leq |\tilde \CP_{m+1}|_{0,\tilde r_{m}}\lesssim \varepsilon_m^{\frac43}.
\end{equation}

With Lemma \ref{lemma_for_conj}, we see that there exists $y_{m+1}\in \Ga^\om$ and $z_{m+1}:S\to \Ga^\om$, satisfying
$$\|y_{m+1}\|_{1,M}<\zeta_{m+1},\qquad \|z_{m+1}\|_{S,1,M}=\left(\sum_{1\leq l\leq k}\|z_{m+1}(\gamma_l)\|^{2}_{1,M}\right)^{\frac12}< \varepsilon_m^{2},$$
such that for every $\gamma\in S$,
\begin{equation}\label{conj_form_m+1}
 \Exp\{y_{m+1}\}\circ \Exp\{P_{m+1}(\gamma)\}\circ \pi(\gamma)\circ \Exp\{y_{m+1}\}^{-1}=\Exp\{z_{m+1}(\gamma)\}\circ \pi(\gamma),
 \end{equation}
 which implies, through Lemma \ref{lem_exp-isometry}, that
 \begin{equation}\label{applem-11}
\Exp\{y_{m+1}\}\circ \Exp\{P_{m+1}(\gamma)\}=\Exp\{z_{m+1}(\gamma)\}\circ \Exp\left\{\pi(\gamma)_* y_{m+1}\right\} .
 \end{equation}
Applying Proposition \ref{lemMoser1} to both sides of (\ref{applem-11}),
we have, on $M$,
\begin{eqnarray*}
& & y_{m+1}+P_{m+1}(\gamma) + s_1(y_{m+1}, P_{m+1}(\gamma)) \\
&=& z_{m+1}(\gamma)+\pi(\gamma)_* y_{m+1} + s_1\left(z_{m+1}(\gamma),\pi(\gamma)_* y_{m+1}\right).
\end{eqnarray*}
By taking $\gamma=\gamma_1,\cdots,\gamma_k$ in the above equality, we obtain that
\begin{equation}\label{eq-homo-refine}
\CH\CP_m + \tilde \CP_{m+1} - d_0y_{m+1}=\CZ_{m+1}+\CQ_{m,1} -\CQ_{m,2},
\end{equation}
where $\CZ_{m+1}$, $\CQ_{m,1}$, $\CQ_{m,2}\in (\Ga^\om)^k$ are defined as
\begin{eqnarray*}
\CZ_{m+1} &:=& \left(z_{m+1}(\gamma_l)\right)_{1\leq l \leq k}, \\
\CQ_{m,1} &:=& \left(s_1 (y_{m+1},P_{m+1}(\gamma_l))\right)_{1\leq l \leq k}, \\
\CQ_{m,2} &:=& \left(s_1(z_{m+1}(\gamma_l),\pi(\gamma_l)_* y_{m+1}\right)_{1\leq l \leq k}.
\end{eqnarray*}
In view of Remark \ref{rmk_L2}, we see that
\begin{equation}\label{small_Zm+1}
\|\CZ_{m+1}\|_{L^2}\lesssim \|z_{m+1}\|_{S,1,M}<\varepsilon_m^{2},
\end{equation}
and, according to Proposition \ref{lemMoser1},
\begin{eqnarray}
& &\|\CQ_{m,1}\|_{L^2} \lesssim \|y_{m+1}\|_{1,M} |\CP_{m+1}|_{0,M} \leq \zeta_{m+1}|\CP_{m+1}|_{0,\tilde r_{m}},\label{small_Qm1}\\
& &\|\CQ_{m,2}\|_{L^2} \lesssim \|z_{m+1}\|_{S,1,M} |y_{m+1}|_{0,M}\lesssim\varepsilon_m^{2}.\label{small_Qm2}
\end{eqnarray}
Through Lemma \ref{lemmaH_d0}, we have  $\CH\circ d_0=0$. Hence, by projecting Eq. (\ref{eq-homo-refine}) onto $\Ker\square$, we obtain
\begin{eqnarray*}
\|\CH\CP_m\|_{L^2}&\leq&\|\CH\tilde \CP_{m+1}\|_{L^2}+  \|\CH\CZ_{m+1}\|_{L^2}+\|\CH\CQ_{m,1}\|_{L^2}+\|\CH\CQ_{m,2}\|_{L^2}\\
&\lesssim&\varepsilon^{\frac43}_m + \zeta_{m+1}|\CP_{m+1}|_{0,\tilde r_{m}}.
\end{eqnarray*}
Then, through Lemma \ref{lem-normharmonic} and the last inequality of (\ref{Pm+1L2}), we have
\begin{eqnarray*}
|\CP_{m+1}|_{0,\tilde r_{m}}&\leq&|\CH\CP_m|_{0,\tilde r_{m}}+ |\tilde \CP_{m+1}|_{0,\tilde r_{m}}\\
&\lesssim&(1+\CJ)e^{\tilde r_{m}\CJ}\|\CH\CP_m\|_{L^2}+ |\tilde \CP_{m+1}|_{0,\tilde r_{m}}\\
&\lesssim&(1+\CJ)e^{\tilde r_{m}\CJ}\left(\varepsilon^{\frac43}_m + \zeta_{m+1}|\CP_{m+1}|_{0,\tilde r_{m}}\right)+\varepsilon^{\frac43}_m,
\end{eqnarray*}
which implies, through (\ref{use_vareps0}), that
\begin{eqnarray*}
|\CP_{m+1}|_{0,\tilde r_{m}}&\leq&|\ln(\varepsilon_0)|^{\frac1{12}} |\ln(\varepsilon_0)|^{\frac16}\left( \varepsilon^{\frac43}_m + \zeta_{m+1} |\CP_{m+1}|_{0,\tilde r_{m}}\right)+|\ln(\varepsilon_0)|^{\frac1{12}} \varepsilon^{\frac43}_m\\
&\leq&2|\ln(\varepsilon_0)|^{\frac14} \varepsilon^{\frac43}_m+|\ln(\varepsilon_0)|^{\frac14} \zeta_{m+1} |\CP_{m+1}|_{0,\tilde r_{m}}
\end{eqnarray*}

Then, by (\ref{zeta01}), $|\ln(\varepsilon_0)|^{\frac14}\zeta_{m+1}<\frac14$. Therefore, (\ref{Pm+10r}) is shown by
$$|\CP_{m+1}|_{0,r''_{m}} \leq \frac{ 2|\ln(\varepsilon_0)|^{\frac14}\varepsilon^{\frac43}_m}{1- |\ln(\varepsilon_0)|^{\frac14}\zeta_{m+1}}\leq \varepsilon^{\frac54}_m.\qed$$

\subsection{Convergence}

The proof of Theorem \ref{thm_conj} and \ref{thm-geom}  is completed by the following convergence argument with the sequences $\{w_m\}$, $\{P_m\}$ constructed as above, satisfying
$$\|w_m\|_{r_m}\leq  \varepsilon_m^{\frac{19}{20}}, \quad \|P_m\|_{S,r_m}\leq  \varepsilon_m.$$

\begin{prop} There exists $ W\in  H_{\frac{r_0}2}^2$ with $\|W\|_{\frac{r_0}2}<\varepsilon_0^{\frac34}$ such that
\begin{equation}\label{conj-limit}
\Exp\{W\}^{-1} \circ \pi_0(\gamma)\circ \Exp\{W\}= \pi(\gamma),\quad \forall \ \gamma\in G.
\end{equation}
\end{prop}
\proof At first, according to Proposition \ref{lemMoser1}, for $W_1:=w_0+w_1+s_1(w_0,w_1)\in H_{r_2}^2$, we have $\Exp\{w_0\}\circ\Exp\{w_1\}=\Exp\{W_1\}$, which implies that, for every $\gamma\in S$,
\begin{eqnarray*}
& & \Exp\{W_1\}^{-1}\circ \Exp\{P_0(\gamma)\} \circ \pi(\gamma) \circ\Exp\{W_1\}\\
&=&\Exp\{w_1\}^{-1}\circ\left(\Exp\{w_0\}^{-1}\circ \Exp\{P_0(\gamma)\} \circ \pi(\gamma)\circ\Exp\{w_0\} \right)\circ\Exp\{w_1\} \\
&=&\Exp\{w_1\}^{-1}\circ \Exp\{P_1(\gamma)\} \circ \pi(\gamma) \circ \Exp\{w_1\} \\
&=& \Exp\{P_2(\gamma)\} \circ \pi(\gamma).
\end{eqnarray*}
The estimates $\|w_0\|_{r_0}\leq \varepsilon_0^{\frac{19}{20}}$ and $\|w_1\|_{r_1}\leq \varepsilon_1^{\frac{19}{20}}$ imply, through (\ref{supnormequiv-1}) and Lemma \ref{lem_norm-1}, that
\begin{eqnarray*}
|s_1(w_0,w_1)|_{0,\tilde r_1}\lesssim |w_0|_{1,\tilde r_0}|w_1|_{0,\tilde r_1}\lesssim \frac{\varepsilon_0^{\frac{19}{20}}}{(r_0-r_1)^{3n+1}}\cdot\frac{\varepsilon_1^{\frac{19}{20}}}{(r_1-r_2)^{3n}}\leq \varepsilon_0^{\frac{19}{20}}\varepsilon_1^{\frac{9}{10}},
\end{eqnarray*}
since by (\ref{vareps0_small}) and (\ref{seqs}), we have
$$ \frac{1}{(r_0-r_1)^{3n+1}(r_1-r_2)^{3n}}=\frac{2^{\frac{5\cdot 5n}{2} +2}}{r_0^{5n+1}}\leq  \varepsilon_1^{-\frac{1}{20}}.$$
Then, by (\ref{supnormequiv-2}), $\|s_1(w_0,w_1)\|_{r_2}\lesssim \varepsilon_0^{\frac{19}{20}}\varepsilon_1^{\frac{9}{10}}$. Hence, we have
$$\|W_1-w_0\|_{r_{2}}\leq \|w_1\|_{r_{1}}+\|s_1(w_0,w_1)\|_{r_{2}}<\varepsilon_1^{\frac{19}{20}}+ |\ln(\varepsilon_0)|^{\frac14}\varepsilon_0^{\frac{19}{20}}\varepsilon_1^{\frac{9}{10}}<\varepsilon_1^{\frac78}.$$

Let us define the sequence of vector fields $\{W_m\}$ by induction. More precisely, assume that there exists $W_m$ with $\|W_m\|_{r_{m+1}}<\varepsilon_0^{\frac56}$ such that, for every $\gamma\in S$,
$$\Exp\{W_{m}\}^{-1}\circ \Exp\{P_0(\gamma)\} \circ \pi(\gamma) \circ\Exp\{W_{m}\}=\Exp\{P_{m+1}(\gamma)\} \circ \pi(\gamma).$$
With $W_{m+1}:=W_m+w_{m+1}+s_1(W_m, w_{m+1})$, we have, through Proposition \ref{lemMoser1}, that
$\Exp\{W_m\}\circ\Exp\{w_{m+1}\}=\Exp\{W_{m+1}\}$,
which implies that, for every $\gamma\in S$,
\begin{eqnarray*}
& &\Exp\{W_{m+1}\}^{-1}\circ \Exp\{P_0(\gamma)\} \circ \pi(\gamma) \circ\Exp\{W_{m+1}\}\\
&=& \Exp\{w_{m+1}\}^{-1}\circ\left(\Exp\{W_{m}\}^{-1}\circ \Exp\{P_0(\gamma)\} \circ \pi(\gamma)\circ\Exp\{W_{m}\}\right)\circ\Exp\{w_{m+1}\}\\
&=& \Exp\{w_{m+1}\}^{-1}\circ \Exp\{P_{m+1}(\gamma)\} \circ \pi(\gamma) \circ\Exp\{w_{m+1}\}\\
&=&\Exp\{P_{m+2}(\gamma)\} \circ \pi(\gamma).
\end{eqnarray*}
Since $\|w_{m+1}\|_{r_{m+1}}\leq \varepsilon_{m+1}^{\frac{19}{20}}$, we have, through (\ref{supnormequiv-1}) and Lemma \ref{lem_norm-1}, that
\begin{eqnarray*}
|s_1(W_m,w_{m+1})|_{0,\tilde r_{m+1}}&\lesssim& |W_m|_{1,r'_{m+1}} |w_{m+1}|_{0,\tilde r_{m+1}}\\
&\lesssim& \frac{\varepsilon_0^{\frac56}}{(r_{m+1}-r_{m+2})^{3n+1}}\cdot\frac{\varepsilon_{m+1}^{\frac{19}{20}}}{(r_{m+1}-r_{m+2})^{3n}} \ \leq \ \varepsilon_0^{\frac56}\varepsilon_{m+1}^{\frac{9}{10}},
\end{eqnarray*}
where, by (\ref{seqs}) and Lemma \ref{prop_seq},
$$ \frac{1}{(r_{m+1}-r_{m+2})^{3n+1}(r_{m+1}-r_{m+2})^{3n}}\leq \varepsilon_{m+1}^{-\frac{1}{20}}.$$
Then, by (\ref{supnormequiv-2}), $\|s_1(W_m,w_{m+1})\|_{r_{m+2}}\lesssim \varepsilon_0^{\frac56}\varepsilon_{m+1}^{\frac{9}{10}}$. Hence we have
\begin{eqnarray}
\|W_{m+1}-W_m\|_{r_{m+2}}&\leq&\|w_{m+1}\|_{r_{m+1}} + \|s_1(W_m,w_{m+1})\|_{r_{m+2}}\label{error-W-m}\\
&\leq&\varepsilon_{m+1}^{\frac{19}{20}}+|\ln(\varepsilon_0)|\varepsilon_0^{\frac56}\varepsilon_{m+1}^{\frac{9}{10}} \  \leq \ \varepsilon_{m+1}^{\frac34}.\nonumber 
\end{eqnarray}

As $m\to\infty$, we have the convergence of $\{W_m\}$ from (\ref{error-W-m}), and for every $m\in \N^*$,
$$ \|W_{m+1}\|_{r_{m+2}}\leq \|w_0\|_{r_1} + \sum_{j=0}^m \|W_{j+1}-W_j\|_{r_{j+2}}<\varepsilon_0^{\frac{19}{20}}+\sum_{j=0}^m  \varepsilon_{j+1}^{\frac34}<\varepsilon_0^{\frac34}.$$
Hence, for the limite $W:=\lim_{m\to\infty}W_m\in H_{\frac{r_0}{2}}^2$, we have $\|W\|_{\frac{r_0}2}<\varepsilon_0^{\frac34}$. Since
$$\|P_m\|_{S,\frac{r_0}2}\leq \|P_m\|_{S,r_m}<\varepsilon_m\to 0,$$
we have $\Exp\{W\}^{-1} \circ \pi_0(\gamma)\circ \Exp\{W\}= \pi(\gamma)$, for every $\gamma\in S$.
Then (\ref{conj-limit}) is shown since both $\pi$ and $\pi_0$ are $G-$actions.\qed

\appendix

\section{Proof of Proposition \ref{lemMoser1}.}\label{app_proof}

\noindent
{\it Proof of (i).} Let $W=W_i$ be a trivializing coordinate patch as in Section \ref{sec_Moser}. Given $q\in M_{r'}\cap W$, with $z=z(q)$, let $\eta:=\tilde w(z)$, $\xi:=\tilde v(z)$.
If $|\xi|_{\kappa}$ is sufficiently small, then, recalling (\ref{complex-exp}), $\Psi(z,\eta)$ defines (the coordinates of) a point in $M_{r}$ and we apply the composition of flows~:
$$	\zeta=P(z,\xi, \tilde w(\Psi(z,\xi)))
	= \xi+\tilde w(\Psi(z,\xi))+ \theta(z,\xi, \tilde w(\Psi(z,\xi))). $$
Let us set
\begin{equation}\label{s1}
s_1(w,v)(z) := \zeta-\tilde w(z)-\xi
= \left(\tilde w(\Psi(z,\xi))-\tilde w(z)\right)+ \theta(z,\xi, \tilde w(\Psi(z,\xi))).
\end{equation}
	All these quantities are well defined if $|w|_{0,r}$ and $|v|_{0,r}$ are sufficiently small. Let $v_1,v_2$ be two holomorphic small enough vector fields on $W\cap{M}_{r'}$ and let us set
$$\om_1(z):= \tilde w(\Psi(z,\tilde v_1(z))),\quad \om_2(z):= \tilde w(\Psi(z,\tilde v_2(z))).$$
We have that
\begin{eqnarray}
& &	\sup_{z=z(q)\in\Delta_1^n\atop{q\in M_{r'}\cap W}} |\om_1(z)-\om_2(z)|_{\ka}\label{omom'}\\
&\leq& \sup_{\tilde q\in M_{r_*}\cap W}\sup_{\zeta\in \C^n\atop{|\zeta|\leq 1}} |D_{z} \tilde w(z(\tilde q))\zeta|_{\ka}
\sup_{z=z(q)\in  \Delta_1^n\atop{q\in M_{r'}\cap W}}| \Psi(z,\tilde v_1(z))- \Psi(z,\tilde v_2(z))|_{\ka}\nonumber\\
&\lesssim&\sup_{\tilde q\in M_{r_*}\cap W}\sup_{\zeta\in \mathbb{C}^n\atop{|\zeta|\leq 1}}|D_{z} \tilde w(z(\tilde q))\zeta|_{\ka}
\sup_{z=z(q)\in  \Delta_1^n\atop{q\in M_{r'}\cap W}}|\tilde v_1(z)-\tilde v_2(z)|_{\ka}.\nonumber
\end{eqnarray}
On the other hand, we have
\begin{eqnarray*}
(s_1(w,v_1)-s_1(w,v_2))(z) &=& (\tilde w(\Psi(z,\tilde v_1(z)))-\tilde w(z)) - (\tilde w(\Psi(z,\tilde v_2(z)))-\tilde w(z))  \\
   & & + \, \theta(z,\tilde v_1(z), \tilde w(\Psi(z,\tilde v_1(z))))- \theta(z,\tilde v_2(z), \tilde w(\Psi(z,\tilde v_2(z)))).
\end{eqnarray*}
	According to \re{s1}, we have
\begin{eqnarray*}
& & |s_1(w,v_1)-s_1(w,v_2)|_{0,r'}\\
&\leq& |\om_1-\om_2|_{0,r'}+|\theta(z,\tilde v_1(z), \om_1(z))-\theta(z,\tilde v_2(z), \om_2(z))|_{0,r'}\\
&\leq& |\om_1-\om_2|_{0,r'}  +\sup_{t\in[0,1]}|\partial_\xi \theta(z,t\tilde v_1(z)+(1-t)\tilde v_2(z), t\om_1(z)+(1-t)\om_2(z))|_{0,r}|v_1-v_2|_{0,r'} \\
& &+ \, \sup_{t\in[0,1]}|\partial_\eta\theta(z,t\tilde v_1(z)+(1-t)\tilde v_2(z), t\om_1(z)+(1-t)\om_2(z))|_{0,r}|\om_1-\om_2|_{0,r'}.
\end{eqnarray*}
As $|v_1|_{0,r}$, $|v_2|_{0,r}$, $|\om_1|_{0,r}$, $|\om_2|_{0,r}$ are uniformly bounded, it is deduced from \re{omom'} that
$$|s_1(w,v_1)-s_1(w,v_2)|_{0,r'}\lesssim |w|_{1,r}|v_1-v_2|_{0,r'}. \qed$$

\noindent
{\it Proof of (ii).} In view of \cite{Mos69}[Lemma 1], we see the existence of $s_1(w,v)$ with $|s_1(w,v)|_{0,M}\lesssim \|w\|_{1,M} |v|_{0,M}$. It remains to estimate $|s_1(w,v)|_{1,M}$, and it is sufficient to show that, for $x_1,x_2\in M$,
\begin{equation}\label{s1-1M}
|s_1(w, v)\left(x_1\right)-s_1(w, v)\left(x_2\right)| \lesssim \|w\|_{1,M} (1+|v|_{1,M})  |x_1-x_2|.
\end{equation}
According to \cite{Mos69}[P. 437--439], we have
$$s_1(w,v)=w(\Phi(x,v(x)))-w(x)+\theta(x,v(x),w(x)),$$
where $\Phi(x, \xi)=x+\xi+\phi(x, \xi)$ is a $C^{\infty}$-vector function defined for sufficiently small $\xi$ satisfying that
\begin{equation}\label{s1-phi}
\phi(x, 0)=\phi_{\xi}(x, 0)=0,
\end{equation}
and $\theta(x, \xi, \eta)$ is uniquely defined for small $|\xi|,|\eta|$ and in $C^{\infty}$, satisfying that
\beq\label{s1-pho}
\left|\theta_{\xi}\right| \lesssim |\eta|,\quad  \theta(x, 0, \eta)=\theta(x, \xi, 0)=0.
\eeq
%
Then, for $x_1,x_2\in M$,
\begin{align*}
	s_1(w, v)\left(x_1\right)-s_1(w, v)\left(x_2\right)
	= &\left(w\left(\Phi\left(x_1, v\left(x_1\right)\right)\right)-w\left(\Phi\left(x_2, v\left(x_2\right)\right)\right)\right)  + \left(w\left(x_2\right)-w\left(x_1\right)\right) \\
	&  + \, \left(\theta\left(x_1, v\left(x_1\right), w\left(x_1\right)\right)-\theta\left(x_2, v\left(x_2\right), w\left(x_2\right)\right)\right),
\end{align*}
where, the three terms are bounded by
\begin{equation}\label{swv1-term1}
\left|w\left(x_2\right)-w\left(x_1\right)\right|\lesssim |w|_{1,M}  |x_1-x_2|,
\end{equation}
\begin{eqnarray}
& &\left|w\left(\Phi\left(x_1, v\left(x_1\right)\right)\right)-w\left(\Phi\left(x_2, v\left(x_2\right)\right)\right)\right| \label{swv1-term2}\\
&\lesssim& |w|_{1,M} \left|\Phi\left(x_1, v\left(x_1\right)\right)-\Phi\left(x_2, v\left(x_2\right)\right)\right|\nonumber\\
&=&|w|_{1,M}  \left|\left(x_1+v(x_1)+\phi\left(x_1, v\left(x_1\right)\right) \right)- \left(x_2+v(x_2)+\phi\left(x_2, v\left(x_2\right)\right)\right)\right| \nonumber\\
&\leq& |w|_{1,M} \left((1+|v|_{1,M} ) |x_1-x_2| \right.\nonumber\\
& & + \,   \left. \left|\phi\left(x_1, v\left(x_1\right)\right)-\phi\left(x_1, v\left(x_2\right)\right)\right|+\left|\phi\left(x_1, v\left(x_2\right)\right)-\phi\left(x_2, v\left(x_2\right)\right)\right| \right) \nonumber\\
&\lesssim& |w|_{1,M} (1+|v|_{1,M}) |x_1-x_2|,\nonumber
\end{eqnarray}
since (\ref{s1-phi}) implies that
\begin{eqnarray*}
\left|\phi\left(x_1, v\left(x_1\right)\right)-\phi\left(x_1, v\left(x_2\right)\right)\right|&\leq&\sup_{\xi}|\phi_\xi(x_1, \cdot)|_   |v(x_1)-v(x_2)|\\
&\leq&\sup_{x,\xi}|\phi_\xi|\cdot |v|_{1,M}   |x_1-x_2|,\\
\left|\phi\left(x_1, v\left(x_2\right)\right)-\phi\left(x_2, v\left(x_2\right)\right)\right| &\leq&\sup_{x}|\phi_x(\cdot,v\left(x_2\right))|   |x_1-x_2|\\
 &\leq&\sup_{x,\xi}|\phi_{x\xi}|\cdot  |x_1-x_2|,
\end{eqnarray*}
\begin{eqnarray}
& &\left|\theta\left(x_1, v\left(x_1\right), w\left(x_1\right)\right)-\theta\left(x_2, v\left(x_2\right), w\left(x_2\right)\right)\right|\label{swv1-term3}\\
&\leq&\left|\theta\left(x_1, v\left(x_1\right), w\left(x_1\right)\right)-\theta\left(x_2, v\left(x_1\right), w\left(x_1\right)\right)\right|\nonumber\\
& & + \, \left|\theta\left(x_2, v\left(x_1\right), w\left(x_1\right)\right)-\theta\left(x_2, v\left(x_2\right), w\left(x_1\right)\right)\right|\nonumber\\
& & + \, \left|\theta\left(x_2, v\left(x_2\right), w\left(x_1\right)\right)-\theta\left(x_2, v\left(x_2\right), w\left(x_2\right)\right)\right|\nonumber\\
&\lesssim& \|w\|_{1,M} (1+|v|_{1,M}) |x_1-x_2|,\nonumber
\end{eqnarray}
since (\ref{s1-pho}) implies that
\begin{eqnarray*}
\left|\theta\left(x_1, v\left(x_1\right), w\left(x_1\right)\right)-\theta\left(x_2, v\left(x_1\right), w\left(x_1\right)\right)\right|&\leq& \sup_x\left|\theta_x\left(\cdot, v\left(x_1\right), w\left(x_1\right)\right)\right| |x_1-x_2|\\
&\leq& \sup_{x,\eta}\left|\theta_{x\eta}\left(\cdot, v\left(x_1\right), \cdot\right)\right| \left|w\right|_{0,M} |x_1-x_2|\\
&\leq&\sup_{x,\xi,\eta}\left|\theta_{x\eta}\right| \left|w\right|_{0,M} |x_1-x_2|,\\
\left|\theta\left(x_2, v\left(x_1\right), w\left(x_1\right)\right)-\theta\left(x_2, v\left(x_2\right), w\left(x_1\right)\right)\right|&\leq& \sup_\xi\left|\theta_\xi\left(x_2,\cdot, w\left(x_1\right)\right)\right| |v(x_1)-v(x_2)|\\
&\lesssim& \left|w\right|_{0,M}  |v|_{1,M} |x_1-x_2|,\\
\left|\theta\left(x_2, v\left(x_2\right), w\left(x_1\right)\right)-\theta\left(x_2, v\left(x_2\right), w\left(x_2\right)\right)\right|&\leq&\sup_{\eta}| \theta_\eta\left(x_2, v\left(x_2\right), \cdot \right)| |w(x_1)-w(x_2)| \\
&\leq& \sup_{x,\xi,\eta}\left|\theta_{\eta}\right| \left|w\right|_{1,M} |x_1-x_2|.
\end{eqnarray*}
Combining (\ref{swv1-term1}) -- (\ref{swv1-term3}), we obtain (\ref{s1-1M}).\qed


\end{document}